\documentclass[a4paper,12pt,reqno]{article}
\pdfoutput=1 

\usepackage{graphicx}
\usepackage{tabularx}
\usepackage{amsfonts}
\usepackage{amssymb}
\usepackage{pstricks}
\usepackage{amsmath,upgreek}
\usepackage{array}
\usepackage{pgfplots}
\usepackage[
    bookmarksnumbered,
    colorlinks, 
    linkcolor=black, 
    citecolor=black]{hyperref}
    
\usepackage{pstricks}
\usepackage{bm}
\usepackage[top=3cm, bottom=3cm, left=3cm, right=3cm]{geometry}

\pgfplotscreateplotcyclelist{ageplot}{%
{blue,mark=o, mark size={3}},
{red,mark=square},
{black,mark=star},
{brown,mark=triangle},
{green,mark=diamond},}

\tikzset{every mark/.append style={scale=0.7}}
\usepackage{pgfplotstable,booktabs}

\pgfplotstableset{
columns/NDOF2/.style={int detect,column type=r,column name={dof}},
columns/NDOF3/.style={int detect,column type=r,column name={dof}},
columns/NormEx2/.style={int detect,column type=r,
column name=$Q(\textbf{u})$,fixed, zerofill,precision=8,},
columns/NormEx3/.style={int detect,column type=r,
column name=$Q(\textbf{u})$,fixed, zerofill,precision=8,},
columns/Norm2/.style={int detect,column type=r,
column name=$Q(\textbf{u}_{es})$,fixed, zerofill,precision=8,},
columns/Norm3/.style={int detect,column type=r,
column name=$Q(\textbf{u}_{es})$,fixed, zerofill,precision=8,},
columns/Err2/.style={int detect,
column name=$Q(\textbf{e}_{es})$,fixed, sci ,dec sep align, zerofill,precision=6,},
columns/Err3/.style={int detect,
column name=$Q(\textbf{e}_{es})$,fixed,sci ,dec sep align, zerofill,precision=6,},
columns/QOI2/.style={int detect,
column name=QoI,fixed, zerofill,precision=10,},
columns/MOI3/.style={int detect,
column name=QoI,fixed, zerofill,precision=10,},
columns/ErrEx2/.style={int detect,
column name=$Q(\textbf{e})$,fixed,sci ,dec sep align, zerofill,precision=6,},
columns/ErrEx3/.style={int detect,
column name=$Q(\textbf{e})$,fixed,sci ,dec sep align, zerofill,precision=6,},
columns/Efect2/.style={int detect,column name=$\theta$,fixed, zerofill,precision=8,},
columns/Efect3/.style={int detect,
column name=$\theta$,fixed, zerofill,precision=8,},
columns/Efect02/.style={
column name=$\theta_{Q(\textbf{u})}$,fixed, zerofill,precision=8,},
columns/Efect03/.style={int detect,
column name=$\theta_{(Q\textbf{u})}$,fixed, zerofill,precision=8,},
columns/Efect22/.style={int detect,
column name=$\theta_{QoI}$,fixed, zerofill,precision=8,},
columns/Efect23/.style={int detect,
column name=$\theta_{QoI}$,fixed, zerofill,precision=8,},
columns/Convrate2/.style={int detect,
column name=$v_{ex}$,fixed, zerofill,precision=8,},
columns/Convrate3/.style={int detect,
column name=$v_{ex}$,fixed, zerofill,precision=8,},
columns/Convrate2/.style={int detect,
column name=$v_{ex}$,fixed, zerofill,precision=8,},
columns/Convrate3/.style={int detect,
column name=$v_{ex}$,fixed, zerofill,precision=8,},
empty cells with={--}, 
every head row/.style={before row=\toprule,after row=\midrule},
every last row/.style={after row=\bottomrule}
}

\newcommand\vm[1]{\bm{\mathrm{#1}}}         
\newcommand{\norms}[1]{  \lVert #1 \rVert}       
\newcommand{\abs}[1]{\left|#1\right|}       
\newcommand{\dd}{\rm d}                     

\title{Mesh adaptivity driven by goal-oriented locally equilibrated superconvergent patch recovery}
  
\author{O.A. Gonz\'alez-Estrada$^{1}$ \and E. Nadal$\,^{2}$ \and J.J. R\'odenas$^{2}$ \and P. Kerfriden$^{1}$ \and  S.P.A. Bordas$^{*1}$  \and  F.J. Fuenmayor$\,^{2}$}

\begin{document}

\maketitle

\begin{center}\small{
$^{1}$Institute of Mechanics and Advanced Materials (IMAM), Cardiff School of Engineering, Cardiff University, 
Queen's Buildings, The Parade, Cardiff CF24 3AA Wales, UK, e-mail: estradaoag@cardiff.ac.uk (orcid:0000-0002-2778-3389), kerfridenp@cardiff.ac.uk, stephane.bordas@alum.northwestern.edu (orcid:00000-0001-7622-2193)\\
$^{2}$Centro de Investigaci\'on de Tecnolog\'ia de Veh\'iculos (CITV), \\Universitat Polit\`{e}cnica de Val\`{e}ncia, E-46022-Valencia, Spain, \\e-mail: jjrodena@mcm.upv.es, ennaso@upv.es (orcid:0000-0002-2808-298X), ffuenmay@mcm.upv.es}
\end{center}


\begin{abstract}
Goal-oriented error estimates (GOEE) have become popular tools to quantify and control the local error in quantities of interest (QoI), which are often more pertinent than local errors in energy for design purposes (e.g. the mean stress or mean displacement in a particular area, the stress intensity factor for fracture problems). These GOEE are one of the key unsolved problems of advanced engineering applications in, for example, the aerospace industry. This work presents a simple recovery-based error estimation technique for QoIs whose main characteristic is the use of an enhanced version of the Superconvergent Patch Recovery (SPR) technique previously used for error estimation in the energy norm. This enhanced SPR technique is used to recover both the primal and dual solutions. It provides a nearly statically admissible stress field that results in accurate estimations of the local contributions to the discretisation error in the QoI and, therefore, in an accurate estimation of this magnitude. This approach leads to a technique with a reasonable computational cost that could easily be implemented into already available finite element codes, or as an independent postprocessing tool.
\end{abstract}


\section{Introduction}

Assessing the quality of numerical simulations has been a critical area of research for many years. In the context of the finite element method, verification tools have been extensively developed since the late 70s. Historically, these tools have aimed at quantifying a global distance between the finite element solution and the (unknown) exact solution. In the context of linear elasticity, this distance is naturally measured in the energy norm, as the optimal properties underlying the finite element method are established in this sense. Adaptivity schemes have been developed to control this global error measure by local mesh refinement. However, the use of energy estimates for verification is rather limited. Indeed, this particular distance does not necessarily reflect the precision of the actual output of the calculation, for instance, a local stress measure  or the displacement of part of the structure. In the last 15 years, a revolution has happened in this area of research: the development of the so-called error estimates in quantities of interest. The idea is to quantify the effect of local errors on the precision of the output of the calculation, which leads to adaptivity schemes and quality assessment that are more relevant to engineering applications. Technically, errors in quantities of interest rely on the evaluation of the quantity of interest by means of solving an auxiliary ``mechanical'' problem. Once the quantity of interest is globalised, a combination of energy estimates for the initial problem and auxiliary problem can be deployed to estimate the error in the quantity of interest.

The energy estimates developed during the last 40 years can be classified into two categories: explicit estimates that try to link a certain measure of the residual to the error in the displacement field, and estimates that post-process solution fields that are somehow more accurate than the finite element solution. 

The first category of estimates mimics the optimal bounds used to prove the convergence of finite element discretisation schemes. Very easy to implement, the explicit estimates include interpolation constants that are problem dependent and difficult to evaluate. For this reason, they are usually not favoured by engineers. Energy estimates based on post-processing reconstruct an approximation of the exact solution by performing local computations. A first family of such estimates, usually called implicit residual-based approaches, rely on multilevel finite element approaches. A finite element problem with a very fine discretisation is solved over each of the local subdomains (either individual elements \cite{ainsworthoden2000}, patches of elements \cite{babuskarheinboldt1978,cottereaudiez2009} or subdomains consisting of an element and its immediate neighbourhood \cite{larssonhansbo2002}). Depending on how the local problem is linked to the global finite element solution whose accuracy is to be evaluated, different properties can be obtained for the estimate. For instance, the equilibrated element residual approach, the flux-free approach and the constitutive relation error, yield estimates that bound the error from above. Error estimates based on local elements with Dirichlet boundary conditions bound the error from below (see for instance \cite{diezpares2003}). Other implicit residual-based approaches tailor the link between the local and global problems so that the efficiency of the estimate is maximised, without focussing on the bounding properties \cite{larssonhansbo2002}. More detailed discussion about this class of estimates can be found in \cite{verfurth1996, steinramm2003, diezegozcue1998}. 

The second post-processing approach consists of deriving simple smoothing techniques that yield a solution field that converges faster than the finite element solution. In the context of the finite element method, a very popular prototype for such approaches is the Zienkiewicz-Zhu estimate (ZZ)  \cite{zienkiewiczzhu1987}, associated, for example, with nodal averaging or with the superconvergent patch recovery (SPR) \cite{zienkiewiczzhu1992}. In the latter case, the idea is to fit a higher-order polynomial approximation of the stress field over each patch of elements. The distance in the energy norm between this smoothed stress field and the finite element stress field is defined as an error estimate. The success of this approach in the engineering community relies on an intuitive mechanical definition and a certain ease of implementation compared to the first class of estimates, without sacrificing the numerical effectivity. These smoothing techniques were extended to enriched approximations in \cite{bordasduflot2007, bordasduflot2008, rodenasgonzalez2008, duflotbordas2008,prangeloehnert2012,hildlleras2010} and to smoothed finite elements (SFEM) in \cite{gonzaleznatarajan2012}. The role of enrichment and statical admissibility in such procedures was discussed in \cite{gonzalezrodenas2012}.

Some contributions have suggested bridging these two approaches, mostly in an attempt to obtain the guaranteed upper bound that some of the residual-based techniques propose, while retaining the ease of implementation of the ZZ framework. The basic identity is that when the recovered stress field is exactly statically admissible, then the ZZ estimate coincides with the constitutive relation error, and bounds the energy error from above. Following this line of thought, D\'iez \textit{et al.} \cite{diezrodenas2007} and R\'odenas \textit{et al.} \cite{rodenasgonzalez2010} presented a methodology to obtain practical upper bounds of the error in the energy norm  using an SPR-based approach where equilibrium was locally imposed on each patch. Other approaches propose a partial fulfilment of the equilibrium to improve the robustness of the smoothing procedure, which can be seen as a trade-off between exactly equilibrated approaches and more simple empirical averaging. This contribution considers the ZZ-type estimate and the use of one of these smoothing techniques, namely SPR-CX (where C stands for the enforcement of equilibrium constraints and X for the smooth+singular decomposition of the recovered stress field). In this approach, the internal equilibrium and boundary equilibrium are constrained to be locally satisfied in the smoothing-based recovery operation. Such an approach has been shown to yield very good practical global and local effectivities, in particular on the boundaries where basic averaging techniques fail to obtain accurate results. Note that the method is conceptually and practically different from what is proposed by authors working on equilibrated residual (ER) approaches \cite{odenprudhomme2001, steinramm2003, cottereaudiez2009} or constitutive relation error methods \cite{ladeveze2006}. Conceptually, methods using the ER concept focus on upper bounding properties, whilst trying to achieve good effectivities. This is achieved by recovering a globally equilibrated stress field. Here, our focus is effectivity only. The patch-equilibration is performed to ensure that the recovered stress field converges faster than the finite element solution. It is not required to satisfy global equilibrium. Practically, most ER techniques solve for local equilibria using polynomial finite elements, yielding asymptotic error bounds. The use of local divergence-conforming approximation spaces associated with energy minimisation principles \cite{ladevezepelle2005, cottereaudiez2009} permits to achieve strict error bounding. In our case, we aim at achieving superconvergence by fitting our polynomial approximation to the finite element solution at particular points of the patch. The equilibrium is introduced as a constraint to this minimisation problem.

As mentioned previously, the error in quantities of interest can be related to the energy estimates by means of auxiliary problems. This was done successfully in a number of cases, including implicit residual approaches \cite{odenprudhomme2001,ruterstein2006,larssonhansbo2002,ladevezerougeot1999}, dual analysis \cite{almeidapereira2006}  and the ZZ smoothing approaches \cite{cirakramm1998,ruterstein2006} with basic nodal averaging or SPR technique. Note that if an energy estimate with bounding properties is used, then a similar bounding of the error in the quantity of interest can be obtained \cite{odenprudhomme2001,ladevezerougeot1999}. However, it is usually difficult to obtain guaranteed error bounds in local quantities of interest while maintaining the effectivity of the estimate. The necessity of obtaining such bounds in an engineering context is also arguable, 
as the reliability of an a posteriori error estimate, which is quantified by its local effectivity, can be verified beforehand on a number of practical cases.

This contribution proposes to  use  the SPR-CX approach to derive an efficient and simple goal-oriented adaptivity procedure in the context of elasticity for linear quantities of interest. In this context, a simple stress smoothing procedure with partial fulfilment of the equilibrium is applied to both the initial problem and the auxiliary problem associated with the quantity of interest. 
As the enhanced smoothing technique relies on the explicit definition of the load, we  emphasise that the auxiliary or dual problem can often be seen as a mechanical problem with body and boundary loads, and prescribed initial stresses and strains that depend on the particular quantity of interest. We derive these mechanical loads for a number of practical quantities of interest, in relation with the SPR-CX smoothing technique. We also extend this approach to handle singular elasticity problems, where the quantity of interest is a generalised stress intensity factor. The methodology is verified in a number of numerical examples. We also implement the hierarchy of practical upper bounds proposed in \cite{ruterstein2006} with the ZZ/SPR-CX technique as the energy estimate. We show that the respective property that we observe for each of the estimates follows the trend observed by the authors of that work.

The paper is organised as follows. In Section 2 we define the model problem and the finite element approximation used. Section 3 is related to the representation of the error in the energy norm and it is used to introduce the nearly equilibrated recovery technique for obtaining enhanced stress fields. In Section 4 we provide the basics behind error estimation in quantities of interest. We define the dual problem and present the smoothing-based error estimate, in particular, the expressions for the dual loads required for the stress recovery procedure. We present the smoothing-based error estimate  and describe the expressions for the dual loads, which are required for the stress recovery process. In Section 5 we present some numerical examples and, finally, conclusions are drawn in Section 6.

\section{Linear elasticity problem solved by the finite element method}

\subsection{Problem statement}

\label{sec:ProbStatement}

In this section, we introduce the 2D linear elasticity  problem to which the proposed methodology is dedicated. Using Voigt notation, we denote by $ \vm{\sigma}=\{\sigma_x, \sigma_y,\sigma_{xy}\}^T$ the Cauchy stress, by $\vm{u}$ the displacement, and by $\vm{\varepsilon}$ the strain, all these fields being defined over domain $\Omega \subset \mathbb{R}^{2}$, of boundary denoted by $\partial \Omega$. Prescribed tractions, denoted by $\vm{t}$, are imposed over part 
$\Gamma _{N}$ of the boundary, while displacements denoted by $\bar{\vm{u}}$ are prescribed over the complementary part $\Gamma_D$ of the boundary. 
$\vm{b}$ denotes the body load. The elasticity problem takes the following form. We look for $(\vm{\sigma},\vm{u})$ satisfying:
 
\begin{align}
\intertext{\textbullet\  statical admissibility}
  \vm{L}^T \vm{\sigma}  + \vm{b} &=  0    &&  {\rm in }\; \Omega   \label{Eq:IntEq} \\
  \vm{G} \vm{\sigma}  &= \vm{t}     &&  {\rm on }\; \Gamma _{N}         \label{Eq:Neumann}\\
\intertext{where $\vm{L}$ is the differential operator for linear elasticity, and $\vm{G}$ is the projection operator that projects the stress field into tractions over the boundary considering the unit normal $\vm{n}$ to $\Gamma_N$ such that}
 \vm{L}^T = 
\begin{bmatrix}
  \partial/\partial x & 0                   & \partial/\partial y \\ 
  0                   & \partial/\partial y & \partial/\partial x
\end{bmatrix},&
\qquad
 \vm{G} = 
\begin{bmatrix}
  n_x & 0 & n_y \\ 0 & n_y & n_x
\end{bmatrix}
\nonumber
\intertext{\textbullet\ kinematic admissibility}  
  \vm{u}  &= \bar{\vm{u}} &&  {\rm on }\; \Gamma _{D}, \label{Eq:Dirichlet} \\
\intertext{\textbullet\ constitutive relation}  
  \vm{\varepsilon} (\vm{u}) &= \vm{L} \vm{u}  &&  {\rm in }\; \Omega   \\
  \vm{\sigma} &= \vm{D} ( \vm{\varepsilon}(\vm{u}) - \vm{\varepsilon}_0) + \vm{\sigma}_0   &&  {\rm in }\; \Omega \label{Eq:ConstitutiveRel}
\intertext{where operator $\vm{D}$ contains the elasticity coefficients of the usual linear isotropic constitutive law relating stress and strain, and $\vm{\sigma}_0$ and $\vm{\varepsilon}_0$ are initial stresses and strains respectively.} \nonumber
\end{align}
Using the notations introduced in \cite{verfurth1999} the problem above takes the primal variational form: 
\begin{equation}\label{Eq:VariationalForm}
\begin{aligned} 
&\text{Find } \vm{u} \in  (V + \{ \vm{w} \}  )  \text{ such that } \forall \vm{v} \in V:\\
&\int _{\Omega}  \vm{\varepsilon}(\vm{u})^{T}  \vm{D} \vm{\varepsilon}(\vm{v}) \dd \Omega = \\
&\int _{\Omega} \vm{v}^{T} \vm{b}   \dd \Omega + 
\int _{\Gamma_N} \vm{v}^{T} \vm{t}  \dd\Gamma + 
\int _{\Omega} \vm{\varepsilon}(\vm{v})^{T}  \vm{D}  \vm{\varepsilon}_0   \dd \Omega -
\int _{\Omega} \vm{\varepsilon}(\vm{v})^{T} \vm{\sigma}_0    \dd \Omega
\end{aligned} 
\end{equation}
where  $V = \{\vm{v} \;|\; \vm{v} \in  [H^1(\Omega)]^2 , \vm{v}|_{\Gamma_D} = \vm{0} \}$ and $ \vm{w}$ is a particular displacement field satisfying the Dirichlet boundary conditions. 

 \subsection{Finite element discretisation}
  
Let us introduce a classical finite element discretisation scheme for the elasticity problem. The approximate displacement field $\vm{u}^{h}$ is searched for in a space of finite dimension $(V^{h} + \{ \vm{w} \}) \subset (V + \{ \vm{w} \})$ such that $V^{h}$  is spanned by locally supported finite element shape functions.
 
Using the Galerkin framework, the primal variational formulation \eqref{Eq:VariationalForm}  is recast in the form:
\begin{equation}\label{Eq:WeakForm}
\begin{aligned} 
&\text{Find } \vm{u}^h \in  \left(V^h + \{ \vm{w} \} \right)  \text{ such that } \forall \vm{v} \in V^{h}:\\
&\int _{\Omega}  \vm{\varepsilon} (\vm{u}^h)^{T} \vm{D}   \vm{\varepsilon}(\vm{v}) \dd \Omega = \\ 
&\int _{\Omega} \vm{v}^{T} \vm{b}   \dd \Omega + \int _{\Gamma_N} \vm{v}^{T} \vm{t}  \dd\Gamma + \int _{\Omega}    \vm{\varepsilon}(\vm{v})^T \vm{D} \vm{\varepsilon}_0   \dd \Omega -
\int _{\Omega} \vm{\varepsilon}(\vm{v})^{T} \vm{\sigma}_0    \dd \Omega
\end{aligned} 
\end{equation}
which can be solved using classical finite element technology \cite{zienkiewicztaylor2000}.

\section{Smoothing-based error estimates in the energy norm}

\subsection{Zienkiewicz--Zhu estimate}

In the absence of other types of errors, the finite element discretisation error is defined by $\vm{e} = \vm{u}-\vm{u}^h$. To quantify the quality of $\vm{u}^h$ one of the standard approaches is to evaluate the energy norm of $\vm{e}$, which is defined by:
\begin{equation} \label{Eq:errorNorm} 
\norms{\vm{e}}^{2} = \int _{\Omega} \vm{\varepsilon}(\vm{e})^{T} \vm{D}  \vm{\varepsilon}(\vm{e})  \dd\Omega .  
\end{equation}
If we introduce the error in the stress field $\vm{e}_\sigma = \vm{\sigma}- \vm{\sigma}^h $, where   $\vm{\sigma}^h = \vm{D}\left( \vm{\varepsilon}(\vm{u}^h)-\vm{\varepsilon}_0\right) + \vm{\sigma}_0$ is the finite element stress field, and make use of the constitutive relation, the previous expression can be written as
\begin{equation}
\norms{\vm{e}}^{2} =  \int _{\Omega}\vm{e}_\sigma ^{T} \vm{D}^{-1} \vm{e}_\sigma  \dd\Omega
 \end{equation}

Following Zienkiewicz-Zhu \cite{zienkiewiczzhu1987}, an estimate $\mathcal{E}^{zz}$ of the exact error measure $\norms{\vm{e}}$ can be obtained by introducing the approximation 
\begin{equation} \label{Eq:ZZ-estimator} 
\norms{\vm{e}}^{2} \approx ({\mathcal{E}^{zz}})^2 = \int _{\Omega}\left( \vm{e}_\sigma^* \right)^{T} \vm{D}^{-1} \left( \vm{e}_\sigma^* \right) \dd\Omega ,  
\end{equation}
where the approximate stress error $\vm{e}_\sigma^*$ is defined by $\vm{e}_\sigma^* = \vm{\sigma}^*- \vm{\sigma}^h$ and $\vm{\sigma}^*$ is the recovered stress field.

If $\vm{\sigma}^*$ converges to the exact solution at a higher rate (superconvergent) than the FE stress solution, then, the ZZ estimate is asymptotically exact, meaning that the approximate error tends towards the exact error with mesh refinement \cite{zienkiewicztaylor2000}. This estimate does not have guaranteed bounding properties unless the recovered field is statically admissible, in which case it coincides (formally) with the Constitutive Relation Error \cite{ladevezeleguillon1983,pledchamoin2011b} and the Equilibrated Residual Approach \cite{ainsworthoden2000}, which are considered difficult to implement. The most appealing feature of the ZZ approach to error estimation is that very simple recovery techniques based on stress smoothing  permit to obtain good effectivities, which explains the popularity of this method in the engineering community. However, the basic smoothing-based recovery techniques, such as the original SPR, which only use the FE results for the recovery, provide good global effectivities but suffer from a local lack of accuracy. For example, the SPR does not use the information along the boundary,  where the imposed tractions are known and the patches can have a reduced number of elements, thus, leading to a well-known loss in accuracy of the recovered stresses along the Neumann boundary. This local lack of accuracy is a major issue in engineering applications, where boundary stresses are typically of interest, and makes such basic estimates difficult to use to drive adaptivity.
To tackle this problem, trade-offs between exactly equilibrated approaches and  smoothing-techniques were developed. The basics of one such method, the SPR-CX technique, are recalled in the next section as they constitute the corner stone to derive the efficient error estimates in quantities of interest presented in this contribution.

\subsection{Enhanced SPR-based stress recovery}
\label{sec:recovery}

\subsubsection{Stress smoothing}

As noted in previous sections, a widely used technique to control the error in the energy norm in the finite element discretisation is the Zienkiewicz-Zhu error estimator shown in (\ref{Eq:ZZ-estimator}), which is based on the recovery (often called smoothing) of an enhanced stress field $\vm{\sigma}^*$. 

To obtain the recovered field $\vm{\sigma}^*$, first, we define the field $\vm{\sigma}^-$ such that we subtract the initial stress and strain from the field $\vm{\sigma}$:
\begin{equation}\label{Eq:subtractIni}
 \vm{\sigma}^- = \vm{\sigma}-\vm{\sigma}_0 + \vm{D}\vm{\varepsilon_0},
\end{equation}
and perform the smoothing on $\vm{\sigma}^-$. Then, the recovered field is
\begin{equation}
\vm{\sigma}^* = (\vm{\sigma}^{-})^* +\vm{\sigma}_0 - \vm{D}\vm{\varepsilon_0},
\end{equation}
where $(\vm{\sigma}^-)^*$ is the smoothed field that corresponds to $\vm{\sigma}^-$.

In the SPR-CX technique, a patch $\mathcal{P}^{(J)}$ is defined as the set of elements connected to a vertex node $J$. On each patch, a polynomial expansion for each one of the components of the recovered stress field is expressed in the form
\begin{equation}
\hat{\sigma}_{k} ^{*} (\vm{x}) = \vm{p}(\vm{x}) \vm{a}_k \quad k=xx,yy,xy
\end{equation}
where $\vm{p}$ represents a polynomial basis and $\vm{a}_k$ are unknown coefficients to be evaluated. Usually, the polynomial basis is chosen equal to the finite element basis for the displacements. The coefficients $\vm{a}_k$ are evaluated using a least squares approximation to the values of FE stresses evaluated at the integration points of the elements within the patch, $\vm{x}_G \in \mathcal{P}^{(J)}$. 

The recovered stress field coupling the three stress components for the 2D case reads:
\begin{equation}
 \hat{\vm{\sigma}}^{*} (\vm{x}) = 
\begin{Bmatrix}
 \hat{\sigma}_{xx}^{*}(\vm{x})\\
 \hat{\sigma}_{yy}^{*}(\vm{x})\\  
 \hat{\sigma}_{xy}^{*}(\vm{x})                                               
\end{Bmatrix} = 
\vm{P}(\vm{x}) \vm{A} =
\begin{bmatrix}
 \vm{p}(\vm{x}) & \vm{0} & \vm{0} \\
 \vm{0} & \vm{p}(\vm{x}) & \vm{0} \\
 \vm{0} & \vm{0} & \vm{p}(\vm{x}) 
\end{bmatrix}
\begin{Bmatrix}
 \vm{a}_{xx}\\
 \vm{a}_{yy}\\  
 \vm{a}_{xy}  
\end{Bmatrix}
\end{equation}

To obtain the coefficients $\vm{A}$ we solve a linear system of equations resulting from the minimisation of the functional 
\begin{equation}\label{Eq:MinFunctional}
 \mathcal{F}^{(J)}(A) = 
  \int_{\mathcal{P}^{(J)}} (\vm{P}  \vm{A} - \vm{\sigma}^{-h} )^2 \dd\Omega
\end{equation}
where $\vm{\sigma}^{-h} = \vm{D}\vm{\varepsilon}(\vm{u}^h)$. 

To obtain a continuous field,  a partition of unity procedure \cite{blackerbelytschko1994} properly weighting the stress interpolation polynomials, obtained from the different patches formed at the vertex nodes of the element containing point $\vm{x}$, is used. The field $\vm{\sigma}^*$ is interpolated using linear shape functions $N^{(J)}$ associated with the $n_v$ vertex nodes such that
\begin{equation}\label{Eq:conjoint_polynomials} 
 \vm{\sigma}^*(\vm{x}) = \sum_{J=1}^{n_v} N^{(J)}(\vm{x})  \hat{\vm{\sigma}}^{*(J)} (\vm{x}) - \vm{D}\vm{\varepsilon}_0(\vm{x}) + \vm{\sigma}_0(\vm{x}).
\end{equation}
Note that in \eqref{Eq:conjoint_polynomials} we add back the contribution of the initial stresses and strains subtracted in \eqref{Eq:subtractIni}.

\subsubsection{Equilibrium conditions}

The accuracy of such estimates depends on the quality of the recovered field. To improve the quality of the recovered fields, numerical results indicate that when recovering singular fields, statical admissibility and suitably chosen enrichment functions improve the accuracy \cite{gonzalezrodenas2012,gonzalezrodenas2012a}\footnote{The use of enrichment to improve recovery based error estimates for enriched approximations was discussed in detail in some of the first papers discussing derivative recovery techniques for enriched finite element approximations, see References \cite{rodenasgonzalez2010,bordasduflot2007,bordasduflot2008,rodenasgonzalez2008,duflotbordas2008,gonzalezrodenas2012}. A very detailed and clear discussion of a wide variety of error estimators and adaptive procedures for discontinuous failure is provided in \cite{pannachetsluys2009}.}. In this work we consider the SPR-CX recovery technique, which is an enhancement of the error estimator introduced in \cite{diezrodenas2007}, to recover the solutions for both the primal and dual problems ( see Section \ref{sec:errorQoI}). The technique incorporates the ideas in \cite{rodenastur2007} to guarantee locally on patches the exact satisfaction of the equilibrium equations, and the extension in \cite{rodenasgonzalez2008} to singular problems.

Constraint equations are introduced via Lagrange multipliers into the linear system used to solve for the coefficients $\vm{A}$ on each patch in order to enforce the satisfaction of the:
\begin{itemize}
    \item Internal equilibrium equations: We define the constraint equation for the internal equilibrium in the patch as:
    \begin{equation}
     \forall \vm{x}_j \in \mathcal{P}^{(J)} \qquad   \vm{L}^T  \hat{\vm{\sigma}}^{*(J)}(\vm{x}_j)  + \vm{L}^T
     (\vm{\sigma}_0(\vm{x}_j) - \vm{D}\vm{\varepsilon}_0(\vm{x}_j)) +  \hat{\vm{b}}(\vm{x}_j)  := \vm{c}^{\rm int}(\vm{x}_j) = 0    
    \end{equation}
    We consider $\hat{\vm{b}}(\vm{x})$ a polynomial least squares fit of degree $p-1$,  $p$ being the degree of the recovered stress field $\hat{\vm{\sigma}}^{*(J)}$, to the actual body forces $\vm{b}(\vm{x})$. We enforce $\vm{c}^{\rm int}(\vm{x}_j)$ at a sufficient number of $j$ non-aligned points ($nie$) to guarantee the exact representation of $\hat{\vm{b}}(\vm{x})$ and solve the linear system. For example, in the case of $\hat{\vm{b}}(\vm{x})$ being represented by a plane, we choose a set of three non-aligned points in $\mathcal{P}^{(J)}$. In 2D, this constraint will add two equations (one per spatial direction) per point to the linear system assembled for the patch, .i.e. $2 nie $ equations.

    \item Boundary equilibrium equations: A point collocation approach is used to impose the satisfaction of a polynomial approximation to the tractions along the Neumann boundary intersecting the patch. The constraint equation reads 
    \begin{equation}
    \forall \vm{x}_j \in \Gamma_N \cap \mathcal{P}^{(J)}\qquad \vm{G}  \hat{\vm{\sigma}}^{*(J)}(\vm{x}_j) + \vm{G} \vm{L}^T
     (\vm{\sigma}_0(\vm{x}_j) - \vm{D}\vm{\varepsilon}_0(\vm{x}_j)) - \vm{t}(\vm{x}_j) := \vm{c}^{\rm ext}(\vm{x}_j)  = 0      
    \end{equation}
    We enforce $\vm{c}^{\rm ext}(\vm{x}_j)$ in $nbe = p+1$ points along the part of the boundary crossing the patch. In the case that more than one boundary intersects the patch, only one curve is considered in order to avoid over-constraining. The boundary constraint adds $2 nbe$ equations to the linear system of equations to solve.  
    
    \item Compatibility equations: $\vm{c}^{\rm cmp}(\vm{x}_j)$  is only imposed in the case that $p \geq 2$ in a sufficient number of non-aligned points $nc$. For example, for $p=2$ we have $nc=1$.  $\hat{\vm{\sigma}}^*$ directly satisfies $\vm{c}^{cmp}$ for $p=1$.
    The  2D compatibility equation   expressed in terms of stresses 
    is:
\begin{equation}
\forall \vm{x}_j \in \mathcal{P}^{(J)} \qquad    
\frac{\partial^2}{\partial y^2} \left( k \hat{\sigma}_{x x}-
\nu q \hat{\sigma}_{y y} \right)+
\frac{\partial^2}{\partial x^2} \left( k \hat{\sigma}_{y y} -
\nu q \hat{\sigma}_{x x} \right)  -
2(1+\nu)\frac{\partial^2 \hat{\sigma}_{x y}} {{\partial x}{\partial y}}=0
\label{eq:Compatibility}
\end{equation}
where $k,q$ are functions of the Poisson's coefficient $\nu$
\[
  \begin{cases}
   k=1, \; q=1 & \text{for plane stress}\\
   k=(1-\nu)^2, \; q=(1+\nu) & \text{for plane strain}
  \end{cases}
\]
Enforcement of compatibility adds $nc$ equations to the linear system for the patch.

\end{itemize}

Thus, after adding the $2nie+nc$ equations and the additional $2nbe$ equations in patches intersected by the Neumann boundary, the Lagrange functional enforcing the constraint equations for a patch $\mathcal{P}^{(J)}$ can be written  as
\begin{equation}
 \mathcal{L}^{(J)} (\vm{A},\vm{\lambda}) = 
 \mathcal{F}^{(J)}(\vm{A}) +
 \sum_{i=1}^{nie}\lambda^{\rm int}_i\left( \vm{c}^{\rm int}(\vm{x}_i)\right) + \sum_{j=1}^{nbe}\lambda^{\rm ext}_j\left(\vm{c}^{\rm ext}(\vm{x}_j)\right) +
\sum_{k=1}^{nc}\lambda^{\rm cmp}_k\left(\vm{c}^{\rm cmp}(\vm{x}_k)\right)
\end{equation}

\begin{figure}[htb!]
    \centering
    \includegraphics{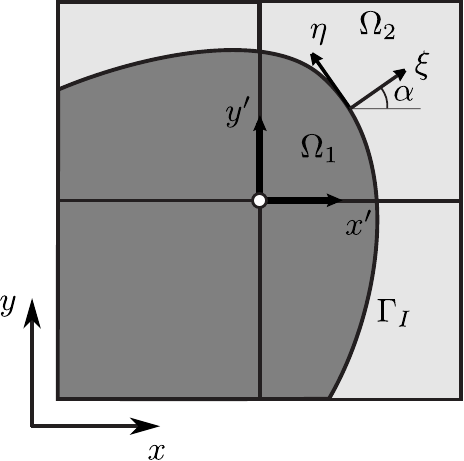}
    \caption{Equilibrium conditions along internal boundaries.}
    \label{fig:InternalBoundary}
\end{figure}

Problems involving an internal interface $\Gamma_I$ along which continuity of the field is not maintained (e.g. bimaterial problems), we use a different polynomial expansion on each side of the boundary and enforce statical admissibility (traction continuity along $\Gamma_I$). Suppose that we have a patch intersected by $\Gamma_I$ such that $\Omega_e = \Omega_{1,e} \cup \Omega_{2,e}$ for intersected elements, as indicated in Figure \ref{fig:InternalBoundary}. To enforce equilibrium conditions along $\Gamma_I$ we define the stresses $\hat{\vm{\sigma}}^*_{\Omega_1}|_{\Gamma_I}$, $\hat{\vm{\sigma}}^*_{\Omega_2}|_{\Gamma_I}$ on each side of the internal boundary. Then, the boundary equilibrium along $\Gamma_I$ given the prescribed tractions $ \vm{t}_{\Gamma_I}= [t_x \; t_y]^T$ is:
\begin{equation}
\vm{G}(\hat{\vm{\sigma}}^*_{\Omega_1}|_{\Gamma_I} -\hat{\vm{\sigma}}^*_{\Omega_2}|_{\Gamma_I})  = \vm{t}_{\Gamma_I}
\label{eq:ContIntEqDispl}
\end{equation}

Once we have the equilibrated recovered fields on each of the patches $\hat{\vm{\sigma}}^{*(J)} $,  \eqref{Eq:conjoint_polynomials} yields a continuous recovered stress field $\vm{\sigma}^{*}$. Note that this postprocess of the local solutions introduces in the global solution $\vm{\sigma}^{*}$ a small lack of equilibrium $\vm{s} = \sum _{J=1} ^{n_v} \nabla N^{(J)} \hat{\vm{\sigma}}^{*(J)}$ when evaluating the divergence of the internal equilibrium equation, as explained in \cite{diezrodenas2007, rodenasgonzalez2010}. 

\subsubsection{Singular fields}

 Different techniques were proposed to account for the eventual non smooth, e.g. singular, part during the recovery process \cite{rodenasgonzalez2008, bordasduflot2007, duflotbordas2008}. Here, following the ideas in \cite{rodenasgonzalez2008} for singular problems, the exact stress field $\vm{\sigma}$  is decomposed into two stress fields, a smooth field $\vm{\sigma}_{\rm smo}$ and a singular field $\vm{\sigma}_{\rm sing}$:
 \begin{equation} \label{Eq:splitting} 
 \vm{\sigma} =   \vm{\sigma}_{\rm smo} + \vm{\sigma}_{\rm sing}.
 \end{equation}
 
 Then, on each patch, the recovered field $\hat{\vm{\sigma}}^*$ required to compute the error estimate given in (\ref{Eq:ZZ-estimator}) can be expressed as the contribution of two recovered stress fields, one smooth $\hat{\vm{\sigma}}^*_{\rm smo}$ and one singular $\hat{\vm{\sigma}}^*_{\rm sing}$:
 \begin{equation} \label{Eq:recovered_splitting} 
 \hat{\vm{\sigma}}^* =   \hat{\vm{\sigma}}^*_{\rm smo} + \hat{\vm{\sigma}}^*_{\rm sing}.
 \end{equation}
 
 For the recovery of the singular part, the expression which describes the asymptotic fields near the crack tip with respect to a coordinate system $(r,\phi)$ at the tip is used \cite{szabobabuska1991}: 
 \begin{equation}\label{Eq:CrackTipField:Stress}
 \vm{\sigma}_{\rm sing}(r,\phi) = 
K_{{\rm I}}  \lambda _{{\rm I}}  r^{\lambda _{{\rm I}} -1}  \vm{\Phi} _{{\rm I}} (\lambda _{{\rm I}} ,\phi ) + 
K_{{\rm II}}  \lambda _{{\rm II}} r^{\lambda _{{\rm II}} -1} \vm{\Phi} _{{\rm II}} (\lambda _{{\rm II}} ,\phi )
 \end{equation}
where $\lambda _m$ (with $m = {\rm I,II}$) are the eigenvalues that determine the order of the singularity, $ \vm{\Phi}_m $ are a set of trigonometric functions that depend on the angular position $\phi$ and $K_m$ are the so-called generalised stress intensity factors (GSIFs)\cite{szabobabuska1991}. To evaluate $\hat{\vm{\sigma}}^*_{\rm sing}$ we first obtain estimated values of the generalised stress intensity factors $K_{\rm I}$ and $K_{\rm II}$ using a domain integral method based on extraction functions \cite{szabobabuska1991, ginertur2009b}. The estimated values of  $K_{\rm I}$ and $K_{\rm II}$ are evaluated once for the singularity and then used to define $\hat{\vm{\sigma}}^*$ on the patches. Notice that the recovered part $\hat{\vm{\sigma}}^*_{\rm sing}$ is an equilibrated field as it satisfies the equilibrium equations.
 
Once the field $\hat{\vm{\sigma}}^*_{\rm sing}$ is evaluated, an FE--type (discontinuous) approximation to the smooth part $\hat{\vm{\sigma}}^h_{\rm smo}$ can be obtained subtracting $\hat{\vm{\sigma}}^*_{\rm sing}$ from the raw FE field:
 \begin{equation} \label{Eq:tens_h_smooth} 
 \hat{\vm{\sigma}}^h_{\rm smo} =   \hat{\vm{\sigma}}^h - \hat{\vm{\sigma}}^*_{\rm sing}.
 \end{equation}
 
 Then, the field $\hat{\vm{\sigma}}^*_{\rm smo}$ is evaluated applying the enhancements of the SPR technique previously described, i.e. satisfaction of equilibrium and compatibility equations on each patch. Note that as both $\hat{\vm{\sigma}}^*_{\rm smo}$  and $\hat{\vm{\sigma}}^*_{\rm sing}$ satisfy the equilibrium equations, $\hat{\vm{\sigma}}^*$ also satisfies equilibrium on each patch. The splitting procedure is only applied on patches within a prescribed distance to the singularity as, away from the singular point,  the solution is smooth enough and can be accurately recovered by the smoothing process.
   
\section{Error estimation in quantities of interest by the SPR-CX recovery technique}
\label{sec:errorQoI}

\subsection{Quantity of interest}

In this section we show how the ZZ estimate may be extended to evaluate the error in a particular quantity of interest \cite{ainsworthoden2000}. For simplicity, in this section we assume that the  Dirichlet boundary conditions are homogeneous. We will consider a quantity of interest $Q: V \rightarrow \mathbb{R}$, defined as a linear functional of the displacement field\footnote{The explanations are restricted to linear quantities of interest. In the developments, affine quantities of the displacement will also be considered, but we will show that this particular case can be recast into the linear case.}. For example, this quantity can be an average of the displacement over a part of the domain. The aim is to estimate the error in functional $Q$, which is expressed by
\begin{equation}
Q (\vm{u}) - Q (\vm{u}^h) = Q (\vm{u} - \vm{u}^h) = Q (\vm{e}).
\end{equation}

\subsection{Auxiliary problem and exact error representation}
\label{sec:DualApproach}

Standard procedures to evaluate $Q (\vm{e})$ introduce the \emph{dual} problem (this terminology comes from the optimal control community) or auxiliary problem
\begin{equation}\label{Eq:dual}
\begin{array}{l}
\displaystyle 
\text{Find } \tilde{\vm{u}} \in V  \text{ such that } \forall \vm{v} \in V,
\\
\displaystyle \int _{\Omega}  \vm{\varepsilon}(\vm{v})^{T}  \vm{D}^{} \vm{\varepsilon} (\tilde{\vm{u}}) \dd \Omega = Q(\vm{v}).
\end{array}
\end{equation}
Problem \eqref{Eq:dual} can be seen as the variational form of an auxiliary mechanical problem. Dual field $\tilde{\vm{u}} \in V$ is a displacement that vanishes over $\Gamma_D$. Test function $\vm{v}$ is a virtual displacement. Field $\tilde{\vm{\sigma}} = \vm{D}(\vm{\varepsilon}(\tilde{\vm{u}}) - \tilde{\vm{\varepsilon}}_0 ) +\tilde{\vm{\sigma}}_0$, where $\tilde{\vm{\sigma}}_0$ and $ \tilde{\vm{\varepsilon}}_0$ are known quantities (initial stress and strain) that will be detailed later on, 
can be interpreted as a mechanical stress field. The left-hand side of \eqref{Eq:dual} is the work of internal forces of the auxiliary mechanical problem. As detailed later on, $Q(\vm{v})$ is the work of an abstract external load for the auxiliary mechanical problem.

The primal and dual problem are solved by the finite element method. Here, we will make use of the same finite element space for both problems. Therefore, we will look for an approximation of $\tilde{\vm{u}} \in V$ using the Galerkin approach:
\begin{equation}\label{Eq:dual3}
\begin{array}{l}
\displaystyle \text{Find } \tilde{\vm{u}}^h \in V^h  \text{ such that }  \forall \vm{v} \in V^h,  \\
\displaystyle \int _{\Omega}  \vm{\varepsilon}(\vm{v})^{T}  \vm{D}^{} \vm{\varepsilon} (\tilde{\vm{u}}^h) \dd \Omega = Q(\vm{v}).
\end{array}
\end{equation}

By choosing $\vm{v}=\vm{e}$ in \eqref{Eq:dual}, and making use of the Galerkin orthogonality for the primal problem (see \cite{ainsworthoden2000} for more details), we can recast expression \eqref{Eq:dual} as follows:
\begin{equation}\label{Eq:dual4}
Q(\vm{e}) = \int _{\Omega}  \vm{\varepsilon}(\vm{e})^{T}  \vm{D}^{} \vm{\varepsilon} (\tilde{\vm{e}}) \dd \Omega 
\end{equation}
where $\tilde{\vm{e}} = \tilde{\vm{u}} - \tilde{\vm{u}}^h$ is the discretisation error of the dual problem \eqref{Eq:dual}. In the context of the ZZ framework, we make use of the constitutive relation in order to obtain an expression in terms of mechanical stresses:
\begin{equation}\label{Eq:errorQoI}
Q(\vm{e}) =  \int_\Omega \vm{e}_{\vm{\sigma}}^{T}  \vm{D}^{-1} \tilde{\vm{e}}_{\vm{\sigma}}  \dd \Omega
\end{equation}
where $\tilde{\vm{e}}_\sigma = \tilde{\vm{\sigma}} - \tilde{\vm{\sigma}}^{h}$ is the stress error of the dual problem and $\tilde{\vm{\sigma}}^{h} = \vm{D}  (\vm{\varepsilon} (\tilde{\vm{u}}^h) - \tilde{\vm{\varepsilon}}_0) + \tilde{\vm{\sigma}}_{0}$ the finite element stress field. 

\subsection{Smoothing-based error estimate}

In order to obtain a computable estimate of the error in the quantity of interest using the ZZ methodology, the stress error fields in expression \eqref{Eq:errorQoI} are replaced by their smoothing-based approximations:
\begin{equation}\label{Eq:errorEstQoI}
Q(\vm{e}) \approx \mathcal{E}_1 = \int_\Omega {\vm{e}_{\vm{\sigma}}^*}^{T}  \vm{D}^{-1} {\tilde{\vm{e}}_{\vm{\sigma}}^*}  \dd \Omega
\end{equation}
where the approximate dual error is $\tilde{\vm{e}}_\sigma^* = \tilde{\vm{\sigma}}^{*} - \tilde{\vm{\sigma}}^{h}$ and $\tilde{\vm{\sigma}}^*$ is the recovered auxiliary stress field. 

The recovered auxiliary stress field can be computed in different ways, for instance by using the standard SPR approach, as proposed in   \cite{cirakramm1998,ruterstein2006}. Here, we propose to recover both the primal and dual stress field by the SPR-CX technique described in Section \ref{sec:recovery}.
Two remarks are worth being made here. First, the analytical expressions that define the tractions and body forces for the dual problem are obtained from the interpretation of the functional $Q$ in terms of tractions, body loads, initial stresses and strains, as seen in Section \ref{sec:DualProblem}. Second, to enforce equilibrium conditions during the recovery process along the boundary of the domain of interest (DoI), we consider it as an internal interface. We use different polynomial expansions on each side of the boundary and enforce statical admissibility of the normal and tangential stresses as previously explained in Section \ref{sec:recovery}.

\subsection{Analytical definitions in the dual problem for enforcing the partial equilibrium}
\label{sec:DualProblem}

The SPR-CX  procedure 
relies on the enforcement of equilibrium for the recovered stress fields evaluated on each patch. For this technique to be applied to recover the dual stress field, the corresponding mechanical equilibrium must be made explicit. In order to do so, the right-hand side of \eqref{Eq:dual} is interpreted as the work of mechanical external forces, and the analytical expression of these forces is derived, depending on the quantity of interest: 
\begin{equation} \label{Eq:WeakFormDual}
\begin{aligned} 
&\text{Find } \tilde{\vm{u}} \in  V  \text{ such that }\forall \vm{v} \in V: \\
&\int _{\Omega}  \vm{\varepsilon} (\vm{v})^{T} \vm{D}   \vm{\varepsilon}(\tilde{\vm{u}}) \dd \Omega = Q(\vm{v})\\ 
&=\int _{\Omega} \vm{v}^{T} \tilde{\vm{b}}   \dd \Omega 
+\int _{\Gamma_N} \vm{v}^{T} \tilde{\vm{t}}  \dd\Gamma 
+\int _{\Omega}    \vm{\varepsilon}(\vm{v})^T \vm{D} \tilde{\vm{\varepsilon}}_0   \dd \Omega 
-\int _{\Omega} \vm{\varepsilon}(\vm{v})^{T} \tilde{\vm{\sigma}}_0    \dd \Omega
\end{aligned} 
\end{equation}

The problem in \eqref{Eq:WeakFormDual} is solved using a FE approximation with test and trial functions in $V^h$. The finite element solution is denoted by $\tilde{\vm{u}}^h \in V^h$.

Such derivations were presented in \cite{rodenas2005,gonzalezrodenas2011b,verdugodiez2011}. Here, we only recall some of the results presented in these papers. Additionally, we provide the analytical expression of the dual load when the quantity of interest is the generalised stress intensity factor (GSIF).

\subsubsection{Mean displacement in \texorpdfstring{$\Omega _I$}{Omegai}}
%
Let us assume that the objective is to evaluate a weighted average of the displacement 
in a sub-domain of interest $\Omega_I \subset \Omega$. In this case, the quantity of interest is
\begin{equation}\label{Eq:umean}
    Q(\vm{u}) 
    = \frac{1}{\left| \Omega_I \right|} \int_{\Omega_I} \vm{u}^T\vm{c}_{u} \dd \Omega
\end{equation}

\noindent where $\left| \Omega_I \right|$ is the measure of $\Omega_I$ and $\vm{c}_{u}$ is an extraction vector used to select the combination of the components of the displacement field that are of interest. For example, $\vm{c}_{u}=\left\{1,0 \right\} ^T$ to extract the first component of $\vm{u}$. 

Now, the right-hand side of expression \eqref{Eq:dual} becomes
\begin{equation}
    Q(\vm{v}) =  
    \int_{\Omega_I} \vm{v}^T \left(  \frac{\vm{c}_{u}}{\left| \Omega_I \right|} \right)  \dd \Omega 
\end{equation}

By comparing this expression to the variational form \eqref{Eq:WeakFormDual} of a general mechanical problem, it is clear that quantity $\tilde{\vm{b}}$ defined by $\tilde{\vm{b}} = {\vm{c}_{u}} / {\left| \Omega_I \right|}$ corresponds to a field  of body forces for the auxiliary mechanical problem.

\subsubsection{Mean displacement along \texorpdfstring{$\Gamma _I$}{Gammai}}

The quantity of interest is now the mean value of the displacement along a subset $\Gamma _I$ of the Neumann boundary $\Gamma_N$:
\begin{equation}\label{Eq:umeanGamma}
    Q(\vm{u}) 
    =\frac{1}{\left| \Gamma_I \right|} \int_{\Gamma_I} \vm{u}^T\vm{c}_{u}  \dd \Gamma
\end{equation}

\noindent where $\left| \Gamma_I \right|$ is the length of $\Gamma_I$ and $\vm{c}_{u }$ an extractor acting on $\vm{u}$. 
%
%
Note that the quantity $\tilde{\vm{t}} = {\vm{c}_{u }} / {\left| \Gamma_I \right|}$ can be interpreted as a vector of tractions applied along the boundary in the problem defined in (\ref{Eq:WeakFormDual}). Thus,  $\tilde{\vm{t}}$ is a vector of tractions applied on $\Gamma_I$ which can be used in the dual problem to extract the mean displacements along $\Gamma_I$.

\subsubsection{Mean strain in \texorpdfstring{$\Omega_I$}{Omegai}}

In this case we are interested in some combination of the components of the strain over a subdomain $\Omega_I$ such that the QoI is given by:

\begin{equation}\label{Eq:epsmean}
 Q(\vm{u}) 
  = \frac{1}{|\Omega_I|}\int_{\Omega_I} \vm{c}_{\varepsilon}^T \vm{\varepsilon}(\vm{u}) \dd \Omega  
  = \int_{\Omega_I} \frac{\vm{c}_{\varepsilon}^T}{|\Omega_I|} \vm{\varepsilon}(\vm{u}) \dd \Omega
\end{equation}
where $\vm{c}_{\varepsilon}$ is the extraction operator used to define the combination of strains under consideration. Thus, the term $\tilde{\vm{\sigma}}_0 = \vm{c}_{\varepsilon}^T/|\Omega_I|$ represents the vector of initial stresses that are used to define the auxiliary problem for this particular QoI.

\subsubsection{Mean stress value in \texorpdfstring{$\Omega_I$}{Omegai}}
\label{sec:MeanStressOmega}

Let us consider now as QoI the mean stress value 
given by a combination $\vm{c}_{\sigma}$ of the stress components $
\vm{\sigma} = \vm{D} ( \vm{\varepsilon}(\vm{u}) - \vm{\varepsilon}_0) + \vm{\sigma}_0$ in a domain of interest which reads:
\begin{equation}
  Q(\textbf{u}) 
  = \frac{1}{|\Omega_I|}\int_{\Omega_I} \vm{c}^T_{\sigma} (\vm{D} ( \vm{\varepsilon}(\vm{u}) - \vm{\varepsilon}_0) + \vm{\sigma}_0)\dd \Omega .
  \label{eq:MoI_mStress}
\end{equation}

$Q$ is an affine functional. Let us define
\begin{equation}
 \tilde{Q}(\vm{v}) = \int _{\Omega}   \vm{c}^T_{\sigma} \vm{D} ( \vm{\varepsilon}(\vm{v}))\dd \Omega  
\end{equation}
for $\vm{v}$ an arbitrary vector of $H^1(\Omega)$. $\tilde{Q}$ is such that $\tilde{Q}(\vm{e}) = Q(\vm{e})$, so that by solving the dual problem
\begin{equation}
 \int _{\Omega}  \vm{\varepsilon} (\vm{v})^{T} \vm{D}   \vm{\varepsilon}(\tilde{\vm{u}}) \dd \Omega = 
\tilde{Q}(\vm{v})
\end{equation} 
for $\tilde{\vm{u}}$, the derivations of Section  \ref{sec:DualApproach} apply.

Similarly to the previous quantity, the right-hand side of the auxiliary problem is defined by the  term $\tilde{\vm{\varepsilon}}_0 = \vm{c}_{\sigma }^T/|\Omega_I|$, which represents in this case a vector of initial strains.

\subsubsection{Mean tractions along \texorpdfstring{$\Gamma_I$}{Gammai} included in \texorpdfstring{$\Gamma_D$}{GammaD}}

Let us assume that we want to evaluate, for example, the mean value of a combination of the tractions $\vm{T}_R$ on a part $\Gamma_I$ of the Dirichlet boundary $\Gamma_D$. The application of the principle of virtual work with  test functions $\vm{v} \in H^1$ that do not necessarily vanish over $\Gamma_D$ gives:
\begin{multline}\label{Eq:ReactioTerm}
\int _{\Gamma_D} \vm{v}^T \vm{T}_R \dd \Gamma = 
\int _{\Omega}  \vm{\varepsilon} (\vm{u})^{T} \vm{D}   \vm{\varepsilon}(\vm{v}) \dd \Omega  
-\int _{\Gamma_N} {\vm{v}}^{T} \vm{t}  \dd\Gamma \\
-\int _{\Omega} {\vm{v}}^{T} \vm{b}  \dd\Omega 
+\int _{\Omega} \vm{\varepsilon}(\vm{v})^{T} \vm{\sigma}_0    \dd \Omega
-\int _{\Omega}    \vm{\varepsilon}(\vm{v})^T \vm{D} \vm{\varepsilon}_0   \dd \Omega
\end{multline}

Extracting the quantity $\int \vm{c}_R ^T \vm{T}_R \dd \Gamma$ , where $\vm{c}_R$ is an extractor defined over $\Gamma_I$, is done by defining the prolongation $\vm{\delta} \in H^1(\Omega)$  of extractor $\vm{c}_R$ such that $\vm{\delta}|_{\Gamma_I} = \vm{c}_R/|\Gamma_I|$. For instance, $\vm{\delta}$ is the finite element field that is null at every node that does not belong to $\Gamma_I$. Substituting $\vm{\delta}$ for $\vm{v}$ in \eqref{Eq:ReactioTerm} yields
\begin{multline}
Q(\vm{u}) = 
\int _{\Omega}  \vm{\varepsilon} (\vm{u})^{T} \vm{D}   \vm{\varepsilon}(\vm{\delta}) \dd \Omega  
-\int _{\Gamma_N} \vm{\delta}^{T} \vm{t}  \dd\Gamma \\
-\int _{\Omega} \vm{\delta}^{T} \vm{b}  \dd\Omega 
+\int _{\Omega} \vm{\varepsilon}(\vm{\delta})^{T} \vm{\sigma}_0    \dd \Omega
-\int _{\Omega}    \vm{\varepsilon}(\vm{\delta})^T \vm{D} \vm{\varepsilon}_0   \dd \Omega 
\end{multline}

$Q$ is an affine functional. We define
\begin{equation}
 \tilde{Q}(\vm{v}) = \int _{\Omega}  \vm{\varepsilon} (\vm{v})^{T} \vm{D}   \vm{\varepsilon}(\vm{\delta}) \dd \Omega  
\end{equation}

$\tilde{Q}$ is such that $\tilde{Q}(\vm{e}) = Q(\vm{e})$, so that by solving the dual problem for $\tilde{\vm{u}}$
\begin{equation}
 \int _{\Omega}  \vm{\varepsilon} (\vm{v})^{T} \vm{D}   \vm{\varepsilon}(\tilde{\vm{u}}) \dd \Omega = 
\tilde{Q}(\vm{v})
\end{equation} 
the derivations of Section  \ref{sec:DualApproach} applies. The dual load is an initial strain $\vm{\varepsilon}(\vm{\delta})$. By recalling the variational form in \eqref{Eq:VariationalForm}  we see that the dual problem is a mechanical problem in $\bar{\vm{u}}=\tilde{\vm{u}}-\vm{\delta}$, where $\bar{\vm{u}}|_{\Gamma_D} = -\vm{\delta}|_{\Gamma_D} = -\vm{c}_R/|\Gamma_I|$, hence, a boundary value problem with Dirichlet boundary conditions $\tilde{\vm{u}} = -\vm{c}_R/|\Gamma_I|$ on $\Gamma_D$.

\subsubsection{Generalised stress intensity factor}

Consider the problem of evaluating the generalised stress intensity factor (GSIF), that characterises the singular solution in problems with reentrant corners and cracks, as the quantity of interest. 

From \cite{szabobabuska1991, rodenasginer2006} we take the expression to evaluate the GSIF, which reads
\begin{equation} \label{Eq:InteractionIntegral} 
    Q(\vm{u}) = K = -\frac{1}{C} \int_{\Omega} \left( \sigma_{jk}^{} u_{k}^{\rm aux} -  \sigma_{jk}^{\rm aux}  u_{k}^{} \right) \frac{\partial{q}}{\partial{x_j}} \dd \Omega
\end{equation}
where 
$ u^{\rm aux}$, $\sigma^{\rm aux}$ are the auxiliary fields used to extract the GSIFs in mode I or mode II and $C$ is a constant that is dependent on the geometry and the loading mode. $q$ is an arbitrary $C^0$ function used to define the extraction zone $\Omega_I$ which must take the value of 1 at the singular point and 0 at the boundaries and $x_j$ refers to the local coordinate system defined at the singularity. 

Rearranging terms in the integral in (\ref{Eq:InteractionIntegral}) we obtain:
\begin{multline} 
   Q(\vm{u}) =  K =  
\int_{\Omega_I}  (\vm{\sigma}^{})^T \left(-\frac{1}{C}\right)
\begin{bmatrix}
u_{1}^{\rm aux} q_{,1} \\
u_{2}^{\rm aux} q_{,2}\\
u_{2}^{\rm aux}q_{,1} + u_{1}^{\rm aux}q_{,2}
\end{bmatrix}  - \\
 (\vm{u}_{}^{})^T \left(-\frac{1}{C}\right)
\begin{bmatrix}
\sigma_{11}^{\rm aux}q_{,1} + \sigma_{21}^{\rm aux} q_{,2} \\
\sigma_{12}^{\rm aux}q_{,1} + \sigma_{22}^{\rm aux} q_{,2}
\end{bmatrix} 
  \dd \Omega
\end{multline}
which can be rewritten as a function of initial strains $\tilde{\vm{\varepsilon}}_0$ and body loads $\tilde{\vm{b}}$:
\begin{equation}
 Q(\vm{u}) =  K =  
\int_{\Omega_I}  \vm{\sigma}_{}^{}(\vm{u})^T \tilde{\vm{\varepsilon}}_0 + 
 (\vm{u}^{})^T \tilde{\vm{b}}    \dd \Omega.
\end{equation}

Thus, if we replace $\vm{u}$ with the vector of arbitrary displacements $\vm{v}$, the quantity of interest can be evaluated from
\begin{equation}
  Q(\vm{v}) =  
\int_{\Omega_I} \vm{\sigma}(\vm{v})^T \tilde{\vm{\varepsilon}}_{0} \dd \Omega + 
\int_{\Omega_I} \vm{v}_{}^{T} \tilde{\vm{b}} \dd \Omega .
\end{equation}

Hence, the initial strains and the body loads per unit volume that must be applied in the dual problem to extract the GSIF are defined as
\begin{equation} \label{Eq:IniStrainSIF}
 \tilde{\vm{\varepsilon}}_{0} = 
-\frac{1}{C}
\begin{bmatrix}
u_1^{\rm aux}q_{,1}\\
u_2^{\rm aux}q_{,2}\\ 
u_2^{\rm aux}q_{,1}+u_1^{\rm aux}q_{,2}
\end{bmatrix}  
\quad,\quad 
\tilde{\vm{b}} = \frac{1}{C} 
\begin{bmatrix}
\sigma_{11}^{\rm aux}q_{,1}+\sigma_{21}^{\rm aux}q_{,2}\\
\sigma_{12}^{\rm aux}q_{,1}+\sigma_{22}^{\rm aux}q_{,2}
\end{bmatrix}
\end{equation}

\subsection{Local contributions}

The ZZ error estimate in quantities of interest given in equation \eqref{Eq:errorEstQoI} can be written in terms of local contributions, which proves useful for adaptivity purposes. 
For a discretisation with $n_e$ elements, element $e$ occupying domain $\Omega_e$ such that $\Omega_e: \Omega= \bigcup_{n_e} \Omega_e$, we can write:
\begin{equation}
 \mathcal{E}_1 = 
 \sum_{n_e} \int _{\Omega_e}  
 {\vm{e}_{\vm{\sigma}}^*} ^{T}  
 \vm{D} ^{-1}
 \tilde{\vm{e}}^*_{\vm{\sigma}}  
 \dd \Omega
\end{equation}

\subsection{Related estimates}

In \cite{ruterstein2006}, the authors introduce a hierarchy of practical bounds, by using the triangle and Cauchy-Schwarz inequalities on the exact error measure \eqref{Eq:dual4}. We recall these estimates as they will be used for benchmarking and verification purposes:
\begin{align}
\mathcal{E}_2 &= 
\sum_{n_e} \left| \int _{\Omega_e}  {\vm{e}_{\vm{\sigma}}^*}^T  \vm{D}^{-1} \tilde{\vm{e}}_{\vm{\sigma}}^*  \dd \Omega \right| \label{Eq:errorEstimates2} \\ 
\mathcal{E}_3 &= 
\sum_{n_e} 
\sqrt{ 
\int _{\Omega_e}  
 {\vm{e}^*_{\vm{\sigma}}}^{T}  
 \vm{D} ^{-1}
 \vm{e}^*_{\vm{\sigma}}  
 \dd \Omega }
 \ \
\sqrt{ 
 \int _{\Omega_e}  
 { \tilde{\vm{e}}^{*T}_{\vm{\sigma}} } 
 \vm{D} ^{-1}
 \tilde{\vm{e}}^*_{\vm{\sigma}}  
 \dd \Omega } \label{Eq:errorEstimates3}\\ 
\mathcal{E}_4 &= 
\sqrt{ 
\int _{\Omega}  
 {\vm{e}^*_{\vm{\sigma}}}^{T}  
 \vm{D} ^{-1}
 \vm{e}^*_{\vm{\sigma}}  
 \dd \Omega }
 \ \
\sqrt{ 
 \int _{\Omega}  
 { \tilde{\vm{e}}^{*T}_{\vm{\sigma}} } 
 \vm{D} ^{-1}
 \tilde{\vm{e}}^*_{\vm{\sigma}}  
 \dd \Omega }. \label{Eq:errorEstimates4}
\end{align}

The properties of these error estimates are discussed in \cite{ruterstein2006}. In particular, $\mathcal{E}_4$ is an upper bound if the recovered primal and dual stresses are exactly equilibrated. Such equilibrated recovered stresses were proposed in the context of ZZ estimates in \cite{diezrodenas2007}. Generally speaking, the authors of \cite{ruterstein2006} observe that the more accurate the estimate, the less likely the practical upper bounding property. 
In \cite{ruterstein2006} the authors used the original SPR technique \cite{zienkiewiczzhu1992}. The following section will show how the use of the SPR-CX improves the accuracy of the error estimator in the quantity of interest both at the local and global level.

\section{Numerical results}

In this section 2D benchmark problems with exact solutions are used to investigate the quality of the proposed technique. The first and second problems have a smooth solution whilst the third is a singular problem. For all models we assume plane strain conditions. The \emph{h}-adaptive FE code used in the numerical examples is based on Cartesian meshes independent of the problem geometry \cite{nadalrodenas2011b,nadalrodenas2013}, with bilinear quadrilateral (Q4) elements\footnote{The interested reader is referred to the recent paper by Moumnassi and colleagues \cite{moumnassibelouettar2011} which discusses recent advances in ``ambient space finite elements'' and proposes hybrid level set/FEMs able to handle sharp corners and edges.}. To represent the domain geometry accurately, the integrals in elements cut by the boundary are restricted to the part of the element within the domain as in \cite{belytschkoparimi2003,moescloirec2003}. The integration procedure in these elements is based on the definition of triangular integration subdomains  within $\Omega_{e}$ and aligning with the geometry of the domain. In elements cut by curved boundaries, these triangular subdomains can have curved boundaries. In that case, the actual geometry is reproduced using a transfinite mapping \cite{gordonhall1973, sevillafernandez-mendez2008}. Dirichlet boundary conditions are imposed using Lagrange multipliers following a procedure similar to that described in \cite{bechetmoes2009} and, more recently,  \cite{moumnassibelouettar2011}\footnote{The interested reader may want to refer to \cite{schweitzer2011,menkbordas2011, hiriyurtuminaro2012,gerstenbergertuminaro2012} regarding preconditioning techniques for systems sharing similar features.}.

To assess the performance of the proposed technique we consider the effectivity index of the error estimator $\theta$ defined as the quotient of the estimated error $\mathcal{E}$ in the quantity of interest over the exact error $Q(\vm{e})$: 
\begin{equation} \label{Eq:Effectivity}  
\theta  =\frac{\mathcal{E}}{Q(\vm{e})}   .
\end{equation}  

We can also represent the effectivity in the QoI, $\theta_{QoI}$, defined as the corrected value of the QoI $Q(\vm{u}^h)$ using the error estimate, divided by the exact value $Q(\vm{u})$:
\begin{equation}
  \theta_{QoI} =\frac{Q(\vm{u}^h)+\mathcal{E}}{ Q(\vm{u})} .
\end{equation}

The relative error in the QoI for the exact and estimated error are
\begin{equation}\label{Eq:RelativeErrors}
 \eta^Q(\vm{e}) = \frac{\abs{Q( \vm{e} )}}{\abs{Q( \vm{u} )}} , \quad
 \eta^Q(\vm{e}_{es}) = \frac{\abs{\mathcal{E}}}{\abs{Q(\vm{u}^h)+\mathcal{E}}} .
\end{equation}

The distribution of the local effectivity index $D$ is analysed at the element level, following the definitions described in \cite{rodenastur2007} for the error in the energy norm, adapted here to the case of the error in QoI:

\begin{equation} \label{Eq:LocalEffectivity}  
\begin{array}{ccc} 
{D=\theta ^{e} -1} & {\rm if} & {\theta ^{e} \ge 1} \\ 
{D=1-\dfrac{1}{\theta ^{e} }} & {\rm if} & {\theta ^{e} < 1} 
\end{array}
\qquad \qquad {\rm with} \qquad 
\theta ^{e} =\dfrac{\mathcal{E}^e}{Q(\vm{e}^e)  }   ,
\end{equation}  

\noindent where superscript $^e$ denotes evaluation at the element level. To evaluate $Q(\vm{e}^e)$  and $\mathcal{E}^e$ we use (\ref{Eq:errorQoI}, \ref{Eq:errorEstQoI}) locally on each element. Nonetheless, we should remark that this is only possible 
for some  problems with analytical solutions as the exact value of $\vm{\sigma}$ is unknown in the vast majority of cases, especially for the dual problem.

Once the error in the QoI is estimated, the local error estimates in each element can be used to perform $h$--adaptive analyses using similar techniques to those already available for the error in the energy norm. The refinement of the mesh using the error estimate as the guiding parameter considers a stopping criterion that checks the value of the global estimated error against a prescribed or desired error. If the estimated error is higher than the desired error then the mesh is refined. Several procedures to perform the refinement are available in the literature. To define the size of the elements in the new mesh we follow the adaptive process described in \cite{ladevezemarin1992, coorevitsladeveze1995, ladevezeleguillon1983} which minimises the number of elements in the new mesh for a given accuracy level. This criterion is equivalent to the traditional approach of equally distributing the error in each element of the new mesh as shown in \cite{libettess1995, fuenmayoroliver1996}. These criteria provide the size of the elements in the new mesh as a function of  (\emph{i}) the ratio of the estimated error in the energy norm in the current mesh to the desired error in the new mesh, and (\emph{ii}) the estimated error in the energy norm on each element, which always takes non negative values. They cannot be directly used in goal-oriented adaptivity because the local contributions to the global error in the QoI, evaluated in each element using (\ref{Eq:errorEstQoI}), can take negative values. Thus, for our implementation using \emph{h}--adaptive routines developed for the error in the energy norm we prepare as input the square root of the absolute value of the error in the QoI -- following the relation from the expressions in (\ref{Eq:errorNorm}, \ref{Eq:errorQoI}) -- and the ratio of the estimated error in the QoI in the current mesh to the desired error in the new mesh. 

The refinement technique provides a distribution of the required new element sizes. These sizes are specified for each element of the current mesh, which will be recursively split into 4 new elements until the sizes of the elements are smaller than the required size. A maximum difference of only one refinement level will be allowed between adjacent elements. In these cases $C^0$ continuity will be enforced by  means of multipoint constraints \cite{abelshepard1979, farhatlacour1998}.

\subsection{Problem 1: Thick-wall cylinder subjected to internal pressure. }

\begin{figure}[htb!]
    \centering
    \includegraphics{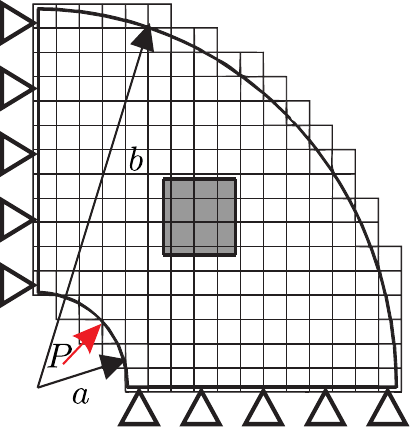}
    \caption{Thick-wall cylinder subjected to an internal pressure. Model and analytical solution, the domain of interest $\Omega_I$ is indicated in dark.}
    \label{fig:Pipe_Model}
\end{figure}

The geometrical model for this problem and the initial mesh are represented in Figure~\ref{fig:Pipe_Model}. Due to symmetry, only $1/4$ of the section is modelled. The domain of interest (DoI) $\Omega_I$ is denoted by a dark square whereas the contours of interest $\Gamma_{\rm i}$ and $\Gamma_{\rm o}$ are the internal and external surfaces of radius $a$ and $b$. Young's modulus is $E=1000 $, Poisson's ratio is $\nu=0.3$, $a=5$, $b=20$ and the load $P=1$.

The exact solution for the radial displacement assuming plane strain conditions is given by
\begin{equation} \label{Eq:uCylinder}  
u_{r}(r) =\frac{P(1+\nu )}{E(c^{2} -1)} \left(r\left(1-2\nu \right) + \frac{b^{2} }{r} \right)  
\end{equation}

\noindent where $c=b/a$, $r=\sqrt{x^{2} +y^{2} }$ and $\phi =\arctan (y/x)$. Stresses in cylindrical coordinates are
\begin{equation} \label{Eq:stressCylinder}  
{\sigma _{r}(r) =\dfrac{P}{c^{2} -1} \left(1-\dfrac{b^{2} }{r^{2} } \right)} \quad
{\sigma _{\phi}(r) =\dfrac{P}{c^{2} -1} \left(1+\dfrac{b^{2} }{r^{2} } \right)} \quad
{\sigma _{z}(r,\phi)  =2 \nu \dfrac{P}{c^2-1}  }
\end{equation} 

Several linear quantities of interest were considered for this problem: the mean radial displacements along $\Gamma_o$, the  mean displacements $u_x$ in the domain of interest $\Omega_I$ and the mean stresses $\sigma_x$ in $\Omega_I$.

\subsubsection{Problem 1.a.: Mean radial  displacements along \texorpdfstring{$\Gamma_o$}{Gammao}}

Let $Q$ be the functional that evaluates the mean normal displacement $\bar{u}_n$ along $\Gamma_o$ such that:
\begin{equation}
Q(\textbf{u}) = \bar{u}_n = 
\frac{1}{\left|\Gamma_o\right|} \int_{\Gamma_o} (\vm{R}\vm{u})^T \vm{c}_u  \dd \Gamma
 \label{eq:MoI_mDisplGamma}
\end{equation}

\noindent where $\vm{R}$ is the standard rotation matrix for the displacements that aligns the coordinate system with the boundary $\Gamma_o$ and $\vm{c}_u= \{1,0\}^T$ is the extraction vector that selects the normal component. The exact value of the QoI given by (\ref{Eq:uCylinder}) for $r=b$ is $\bar{u}_n = 2.42\bar{6} \cdot 10^{-3}$. 

\begin{figure}[htb!]
\centering
\begin{tabular}{c c}
\includegraphics[width=0.3\textwidth]{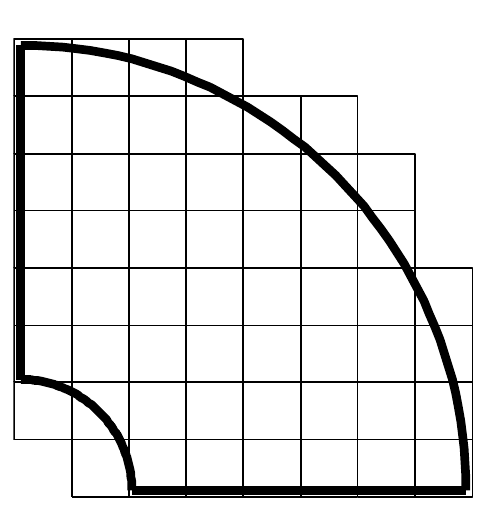} &
\includegraphics[width=0.3\textwidth]{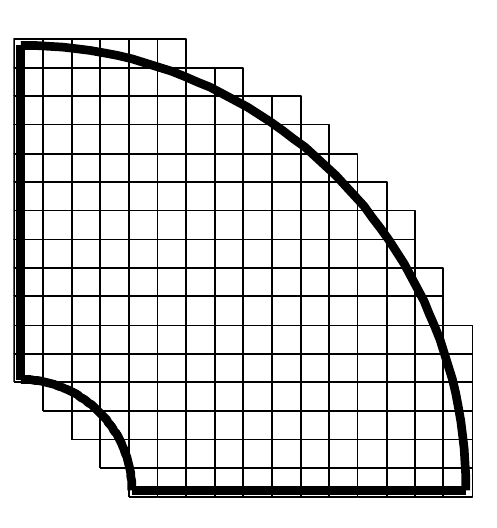}
\\
{a) Mesh 1} & {b)Mesh 2}
\\
\includegraphics[width=0.3\textwidth]{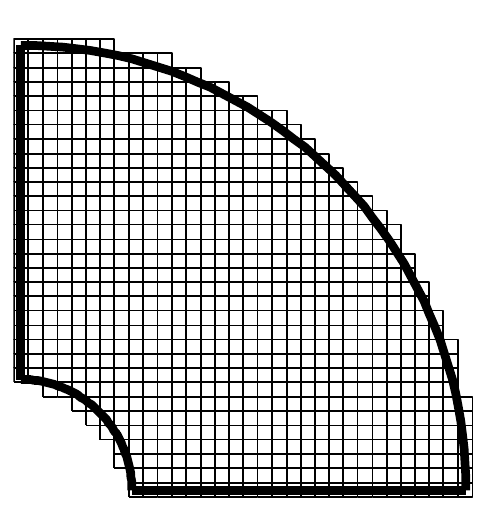} &
\includegraphics[width=0.3\textwidth]{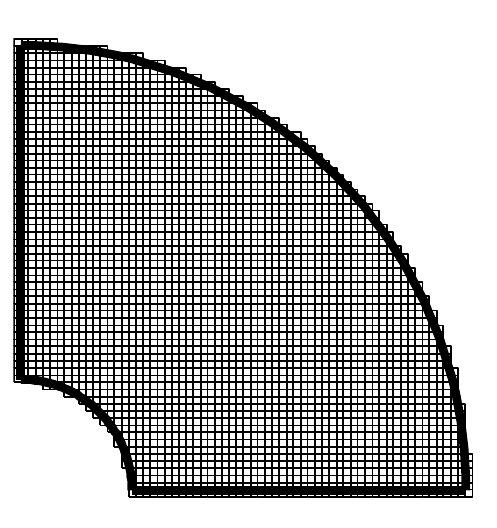}
\\
{c) Mesh 3} & {d)Mesh 4}
\end{tabular}
\caption{Problem 1.a. Sequence of meshes with uniform refinement for the cylinder under internal pressure.}%
\label{fig:CYL_UniformMeshes}%
\end{figure}

\begin{figure}[htb!]
\centering
\begin{tabular}{c c}
\includegraphics[width=0.3\textwidth]{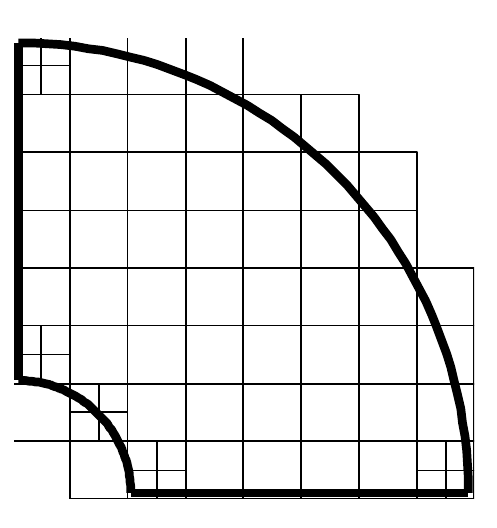} &
\includegraphics[width=0.3\textwidth]{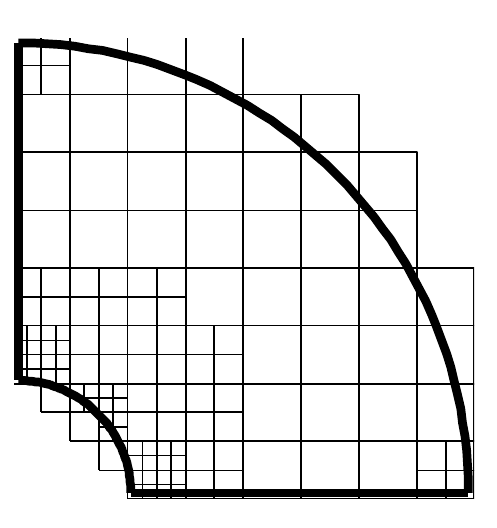}
\\
{a) Mesh 1} & {b)Mesh 2}
\\
\includegraphics[width=0.3\textwidth]{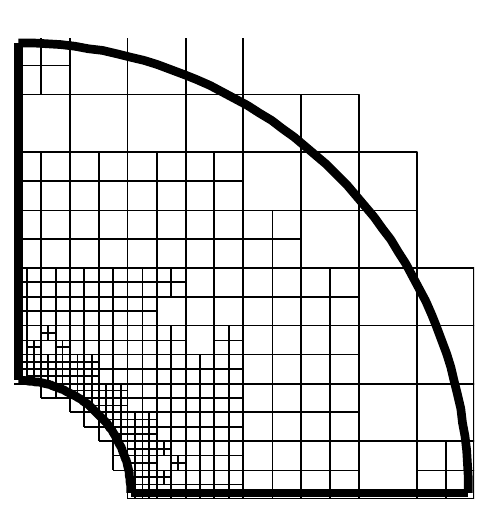} &
\includegraphics[width=0.3\textwidth]{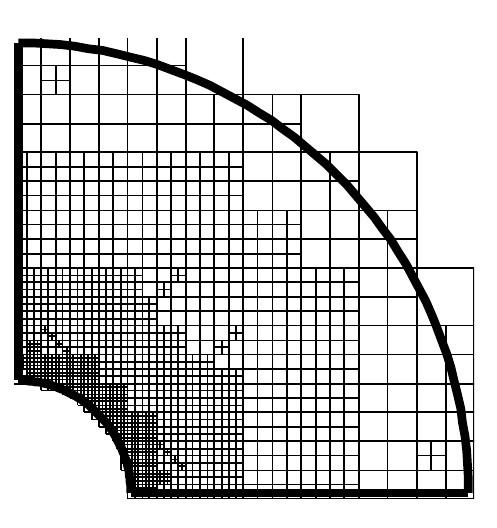}
\\
{c) Mesh 3} & {d)Mesh 4}
\end{tabular}
\caption{Problem 1.a. Sequence of \emph{h}-adapted meshes for the mean radial displacement along $\Gamma_o$.}%
\label{fig:Pipe_mDG_Mesh}%
\end{figure}

To characterise the error before using it in the adaptive process we use the set of uniformly refined meshes shown in Figure \ref{fig:CYL_UniformMeshes}. Figure \ref{fig:Pipe_mDG_Mesh} shows the set of meshes resulting from the \emph{h}-adaptive process guided by the error estimate in this QoI. 

\begin{figure}[htb!]
      \centering
\begin{tikzpicture}
    \begin{semilogxaxis}[
    width=0.48\textwidth,
    title={$\theta$, uniform meshes},
    xtick={100,1000},
    xlabel={DOF},
    legend style={cells={anchor=west}, font=\small},
    legend style={at={(0.98,0.98)},anchor=north east},
    cycle list name=ageplot, table/col sep=comma
    ]a
    \draw (axis cs: 0,1)--(axis cs: 1e8,1); 
    \addplot table[x=dof,y= theta]{Datos/CYLMeanUGammaUniform.csv};
    \addplot table[x=dof,y= theta]{Datos/CYLMeanUGammaSPRUniform.csv};
    \legend{{$\mathcal{E}_1^{\rm SPR-CX}$},{$\mathcal{E}_1^{\rm SPR}$}}
    \end{semilogxaxis}
      \end{tikzpicture} 
      \begin{tikzpicture}
    \begin{semilogxaxis}[
    width=0.48\textwidth,
    title={$\theta$, \emph{h}-adapted meshes},
    xmin=100, xmax=5000,
    xtick={100,1000},
    xlabel={DOF},
    legend style={cells={anchor=west}, font=\small},
    legend style={at={(0.98,0.02)},anchor=south east},
    cycle list name=ageplot, table/col sep=comma
    ]a
    \draw (axis cs: 0,1)--(axis cs: 1e8,1); 
    \addplot table[x=dof,y= theta]{Datos/CYLMeanUGamma.csv};
    \addplot table[x=dof,y= theta]{Datos/CYLMeanUGammaSPR.csv};
    \legend{{$\mathcal{E}_1^{\rm SPR-CX}$},{$\mathcal{E}_1^{\rm SPR}$}}
    \end{semilogxaxis}
      \end{tikzpicture}  
\caption{Problem 1.a. Evolution of the effectivity index $\theta$ considering locally equilibrated , $\mathcal{E}_1^{\rm SPR-CX}$, and non-equilibrated recovery, $\mathcal{E}_1^{\rm SPR}$, for uniform (left) and \emph{h}-adapted (right) meshes.}  %
\label{fig:CylThetaDispGamma}
\end{figure}
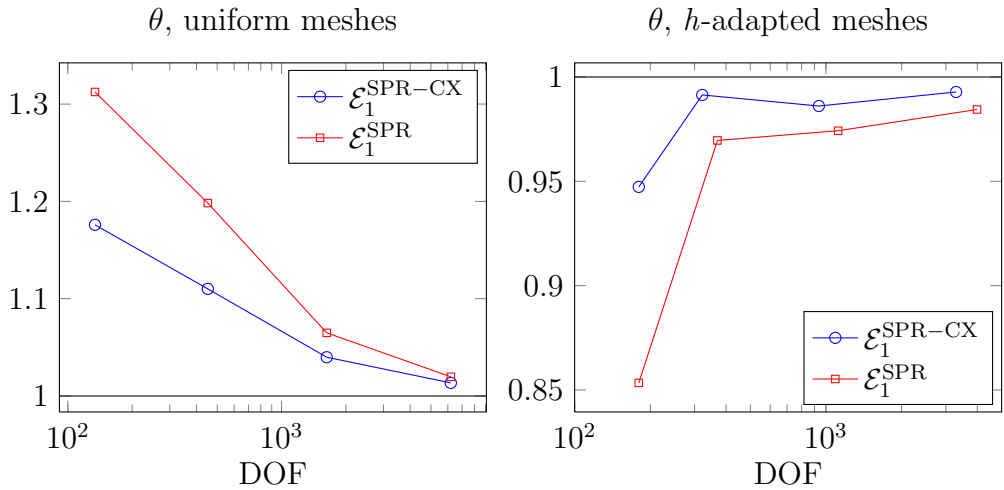

In Tables \ref{tab:MoI_mDG_ConvergQ4Uniform} and \ref{tab:MoI_mDG_ConvergQ4} we show the results for the error estimate $\mathcal{E}_1^{\rm SPR-CX} = \mathcal{E}_1$ from (\ref{Eq:errorEstQoI}), evaluated using the SPR-CX recovery technique  in both the primal and the dual problem, and the exact error $Q(\vm{e})$ for uniform and \emph{h}-adapted meshes, respectively. The recovery technique accurately captures the exact error and provides good effectivities, with values of $\theta = 0.9473$ for the coarsest  \emph{h}-adapted mesh with 180 degrees of freedom (DOF). The evolution of the effectivity index for the proposed technique and the standard SPR is represented in Figure \ref{fig:CylThetaDispGamma}. In this case, SPR-CX converges faster to the optimal value $\theta=1$ and shows a high level of accuracy ($\theta = 0.9928$ for 3294 DOF in \emph{h}-adapted meshes). 

\begin{table}[htb!]
\centering\caption{Problem 1.a. Values for the error estimate $\mathcal{E}_1^{\rm SPR-CX}$ and effectivities considering the uniform meshes in Figure \ref{fig:CYL_UniformMeshes}.}
\begin{filecontents*}{Datos/CYLMeanUGammaUniform.csv}
dof,Qu,Norm,Qees,MOI,Qe,etaes,eta,theta,Efect0,thetaQoI,E2,E3,E4,Eta2,Eta3,Eta4
1.3400000000000e+002,-7.1066666666667e-003,-6.8394450051457e-003,-3.1422686748839e-004,-7.1536718726341e-003,-2.6722166152092e-004,4.4215788108122e-002,3.7601547118328e-002,1.1759034267654e+000,9.6239845288167e-001,1.0066142409898e+000,1.1759034267654e+000,1.1759034267654e+000,1.1759034267654e+000,4.4215788108122e-002,4.4215788108122e-002,4.4215788108122e-002
4.5200000000000e+002,-7.1066666666667e-003,-7.0262964419252e-003,-8.9218182585635e-005,-7.1155146245108e-003,-8.0370224741473e-005,1.2554153271900e-002,1.1309131061183e-002,1.1100899975411e+000,9.8869086893882e-001,1.0012450222107e+000,1.1100899975411e+000,1.1100899975411e+000,1.1100899975411e+000,1.2554153271900e-002,1.2554153271900e-002,1.2554153271900e-002
1.6280000000000e+003,-7.1066666666667e-003,-7.0851356331843e-003,-2.2388471119201e-005,-7.1075241043035e-003,-2.1531033482351e-005,3.1503477184617e-003,3.0296951429200e-003,1.0398233386034e+000,9.9697030485708e-001,1.0001206525755e+000,1.0398233386034e+000,1.0398233386034e+000,1.0398233386034e+000,3.1503477184617e-003,3.1503477184617e-003,3.1503477184617e-003
6.1960000000000e+003,-7.1066666666667e-003,-7.1011091328095e-003,-5.6324364748084e-006,-7.1067415692843e-003,-5.5575338571380e-006,7.9255672722444e-004,7.8201695925956e-004,1.0134776718588e+000,9.9921798304074e-001,1.0000105397680e+000,1.0134776718588e+000,1.0134776718588e+000,1.0134776718588e+000,7.9255672722444e-004,7.9255672722444e-004,7.9255672722444e-004
\end{filecontents*}
\pgfplotstabletypeset[columns={dof,Qees,Qe,theta,thetaQoI},
 columns/dof/.style   ={column name=dof, column type=r, int detect},
 columns/theta/.style ={column name=$\theta$,fixed, zerofill,precision=8},
 columns/thetaQoI/.style ={column name=$\theta_{QoI}$,fixed, zerofill,precision=8},
 columns/Qees/.style  ={column name={$\mathcal{E}_1^{\rm SPR-CX}$},
    fixed,sci ,dec sep align, zerofill,precision=6,},
 columns/Qe/.style    ={column name=$Q(\vm{e})$,
    fixed,sci ,dec sep align, zerofill,precision=6,},col sep= comma
]
{Datos/CYLMeanUGammaUniform.csv} 
\label{tab:MoI_mDG_ConvergQ4Uniform}
\end{table}

\begin{table}[htb!]
\centering\caption{Problem 1.a. Values for the error estimate $\mathcal{E}_1^{\rm SPR-CX}$ and effectivities considering the \textit{h}-adapted meshes in Figure \ref{fig:Pipe_mDG_Mesh}.}
\begin{filecontents*}{Datos/CYLMeanUGamma.csv}
dof,Qu,Norm,Qees,MOI,Qe,etaes,eta,theta,Efect0,thetaQoI,E2,E3,E4,Eta2,Eta3,Eta4,meanD,sigD
1.8000000000000e+002,2.4266666666667e-003,2.3642837534655e-003,5.9096511583934e-005,2.4233802650494e-003,6.2382913201173e-005,2.4352958070302e-002,2.5707244451033e-002,9.4731888190856e-001,9.7429275554897e-001,9.9864571361927e-001,9.4731888190856e-001,9.4731888190856e-001,9.4731888190856e-001,2.4352958070302e-002,2.4352958070302e-002,2.4352958070302e-002,3.6287492971098e-001,5.1859160433818e-001
3.2200000000000e+002,2.4266666666667e-003,2.4091163955624e-003,1.7399131067828e-005,2.4265155266302e-003,1.7550271104290e-005,7.1699715938851e-003,7.2322545759438e-003,9.9138816514205e-001,9.9276774542406e-001,9.9993771701794e-001,9.9138816514205e-001,9.9138816514205e-001,9.9138816514205e-001,7.1699715938851e-003,7.1699715938851e-003,7.1699715938851e-003,1.3201121508096e-001,1.9215889923017e-001
9.3700000000000e+002,2.4266666666667e-003,2.4222073527918e-003,4.3972555262102e-006,2.4266046083180e-003,4.4593138748337e-006,1.8120558487130e-003,1.8376293440249e-003,9.8608343113640e-001,9.9816237065598e-001,9.9997442650469e-001,9.8608343113640e-001,9.8608343113640e-001,9.8608343113640e-001,1.8120558487130e-003,1.8120558487130e-003,1.8120558487130e-003,1.0679105807912e-001,1.7251873679426e-001
3.2940000000000e+003,2.4266666666667e-003,2.4255812555353e-003,1.0776121315049e-006,2.4266588676668e-003,1.0854111313439e-006,4.4407093331246e-004,4.4728480687248e-004,9.9281470438828e-001,9.9955271519313e-001,9.9999678612644e-001,9.9281470438828e-001,9.9281470438828e-001,9.9281470438828e-001,4.4407093331246e-004,4.4407093331246e-004,4.4407093331246e-004,5.1519790264811e-002,1.1414143352401e-001
\end{filecontents*}
\pgfplotstabletypeset[columns={dof,Qees,Qe,theta,thetaQoI},
 columns/dof/.style   ={column name=dof, column type=r, int detect},
 columns/theta/.style ={column name=$\theta$,fixed, zerofill,precision=8},
 columns/thetaQoI/.style ={column name=$\theta_{QoI}$,fixed, zerofill,precision=8},
 columns/Qees/.style  ={column name={$\mathcal{E}_1^{\rm SPR-CX}$},
    fixed,sci ,dec sep align, zerofill,precision=6,},
 columns/Qe/.style    ={column name=$Q(\vm{e})$,
    fixed,sci ,dec sep align, zerofill,precision=6,},col sep= comma
]
{Datos/CYLMeanUGamma.csv} 
\label{tab:MoI_mDG_ConvergQ4}
\end{table}

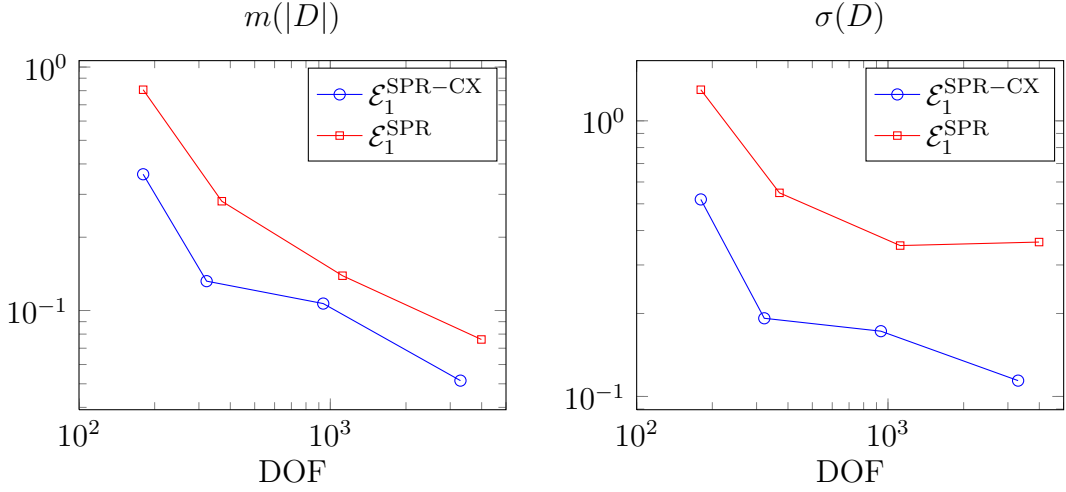
\begin{figure}[htb!]
      \centering
   \begin{minipage}[c]{0.48\textwidth}
       \centering
      \begin{tikzpicture}
    \begin{loglogaxis}[
    width=1\textwidth,
    title={$m(|D|)$},
    xmin=100, xmax=5000,
    xtick={100,1000},
    xlabel={DOF},
    legend style={cells={anchor=west}, font=\small},
    legend style={at={(0.98,0.98)},anchor=north east},
    cycle list name=ageplot, table/col sep=comma
    ]a
    \addplot table[x=dof,y= meanD]{Datos/CYLMeanUGamma.csv};
    \addplot table[x=dof,y= meanD]{Datos/CYLMeanUGammaSPR.csv};
    \legend{{$\mathcal{E}_1^{\rm SPR-CX}$},{$\mathcal{E}_1^{\rm SPR}$}}
    \end{loglogaxis}
      \end{tikzpicture}  
   \end{minipage}
   \begin{minipage}[c]{0.48\textwidth}
       \centering
      \begin{tikzpicture}
    \begin{loglogaxis}[
    width=1\textwidth,
    title={$\sigma(D)$},
    xmin=100, xmax=5000,
    xtick={100,1000},
    xlabel={DOF},
    legend style={cells={anchor=west}, font=\small},
    legend style={at={(0.98,0.98)},anchor=north east},
    cycle list name=ageplot, table/col sep=comma
    ]a
    \addplot table[x=dof,y= sigD]{Datos/CYLMeanUGamma.csv};
    \addplot table[x=dof,y= sigD]{Datos/CYLMeanUGammaSPR.csv};
    \legend{{$\mathcal{E}_1^{\rm SPR-CX}$},{$\mathcal{E}_1^{\rm SPR}$}}
    \end{loglogaxis}
      \end{tikzpicture}  
   \end{minipage}
\caption{Problem 1.a. Evolution of the mean absolute value $m(|D|)$ and standard deviation $\sigma(D)$ of the local effectivity considering locally equilibrated , $\mathcal{E}_1^{\rm SPR-CX}$, and non-equilibrated recovery, $\mathcal{E}_1^{\rm SPR}$.}   
\label{fig:CylmeanDUGamma}
\end{figure}

Note that in this case, the dual problem corresponds to a traction $\tilde{\vm{t}}={\vm{R}^T\vm{c}_u} / {|\Gamma_o|}$ that represents a constant traction normal to the  external boundary. Therefore, the solution to the dual problem can be evaluated exactly for this quantity of interest using the analytical solution for a cylinder under external pressure:
\begin{equation} \label{Eq:stressCylinderOut}  
{\sigma _{r}(r) = -\dfrac{P_o b^2}{b^{2} -a^2} \left(1-\dfrac{a^{2} }{r^{2} } \right)} \quad
{\sigma _{\phi}(r) = -\dfrac{P_o b^2}{b^{2} -a^2} \left(1 + \dfrac{a^{2} }{r^{2} } \right)} 
\end{equation}

\noindent where $P_o$ represents the applied external pressure. As the exact solution for the dual problem is available, it is possible to evaluate the local effectivity $D$ in (\ref{Eq:LocalEffectivity}) at the element level to evaluate the local quality of the error estimate in the quantity of interest. Figure \ref{fig:CylmeanDUGamma} shows the evolution of the mean absolute  value $m(|D|)$ and standard deviation $\sigma(D)$ of the local effectivity for \emph{h}-adapted meshes. The ideal scenario is that both parameters are as small as possible and go to zero as we increase the number of DOFs. In the figure we see that the SPR-CX gives a better local estimation, with values of $\sigma(D)$ and $m(|D|)$ that are smaller than those for the SPR --  for the mesh with approx. 3,300 DOF, $\sigma(D) = 0.11$ for the SPR-CX, compared to $\sigma(D) = 0.36$ for the SPR.

As the \emph{h}-adaptive procedures use local information to refine the mesh, providing an accurate local error estimate to the adaptive algorithm is of critical importance. For this reason, the local performance of the proposed technique indicates that the error estimator based on equilibrated recovered fields for the primal and dual problems is superior to the standard recovery techniques to guide the \emph{h}-adaptive process, even in cases in which the global effectivity is similar to the effectivity of the standard SPR.

\subsubsection{Problem 1.b.: Mean displacements \texorpdfstring{$\bar{u}_x$}{ux} in \texorpdfstring{$\Omega_I$}{Omegai}}

Let us consider the mean displacement $\bar{u}_x$ in $\Omega_I$
as the quantity of interest. The objective is to evaluate the error when evaluating  $\bar{u}_x$ defined by the functional:
\begin{equation}
Q(\textbf{u})=\bar{u}_x=\frac{1}{\left|\Omega_I\right|}\int_{\Omega_I}u_x \dd \Omega
 \label{eq:MoI_mDispl}
\end{equation}

\begin{figure}[htb!]
\centering
\begin{tabular}{c c}
\includegraphics[width=0.3\textwidth]{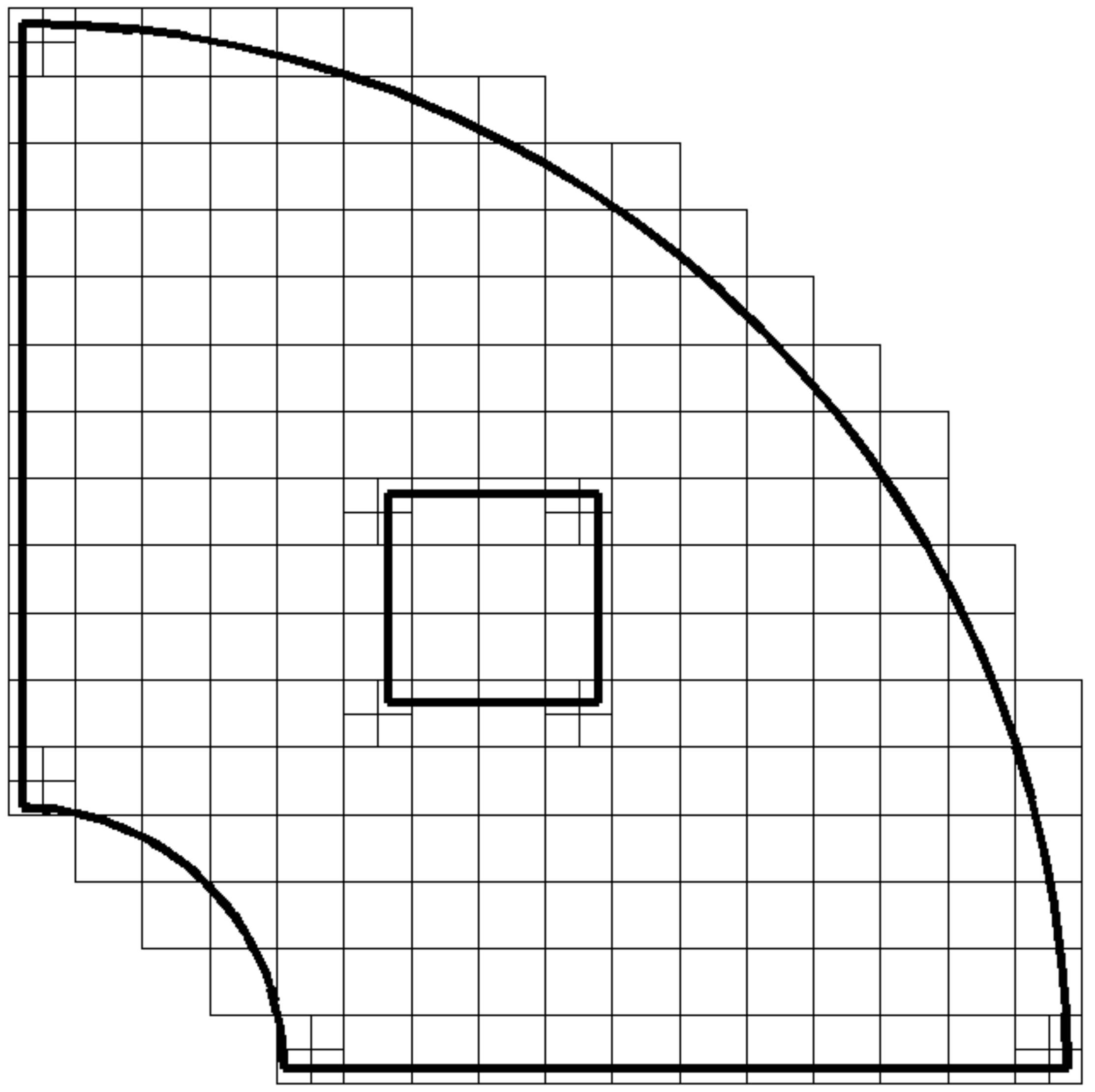} &
\includegraphics[width=0.3\textwidth]{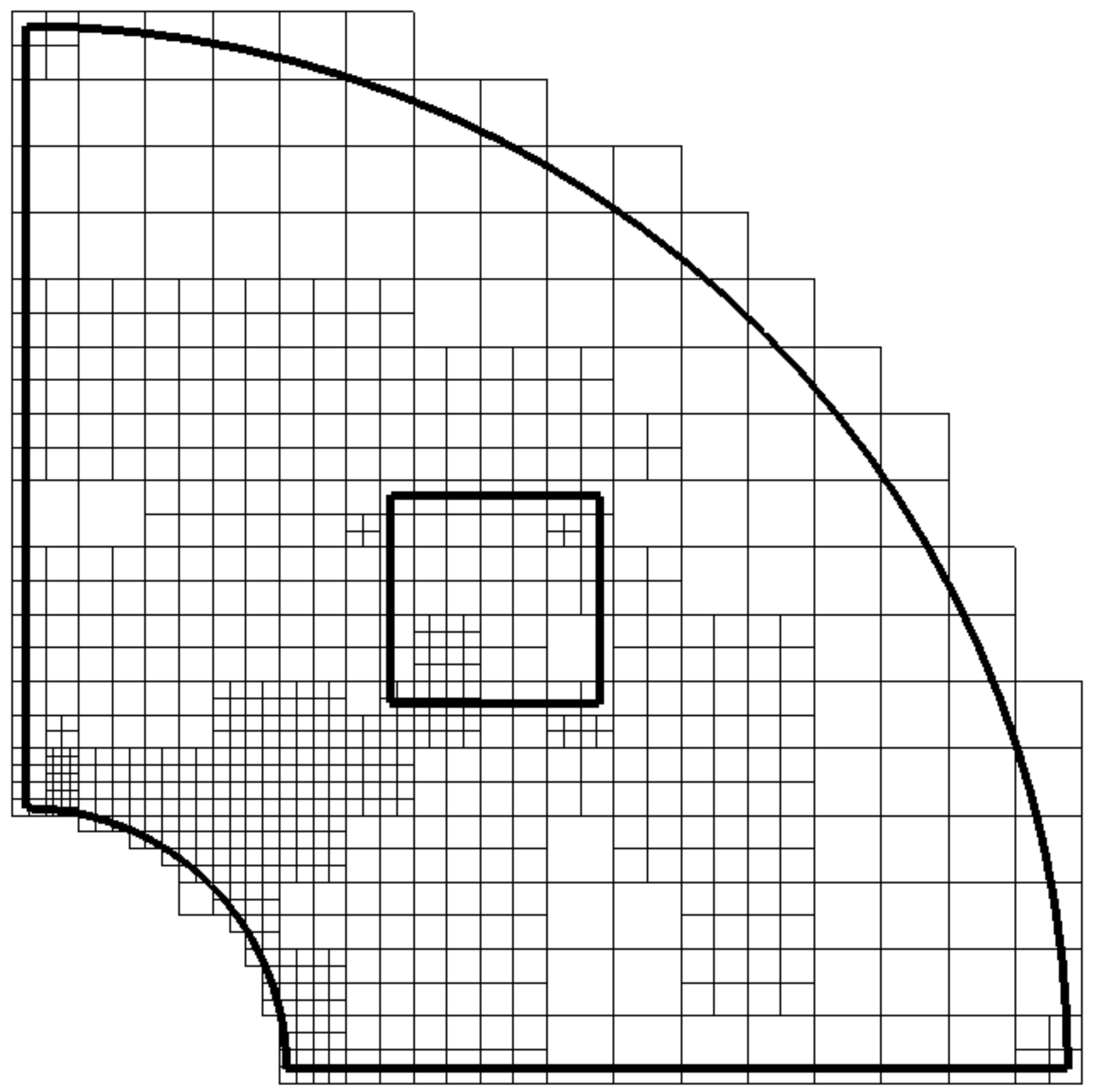}\\
{a) Mesh 1} & {b)Mesh 2}\\
\includegraphics[width=0.3\textwidth]{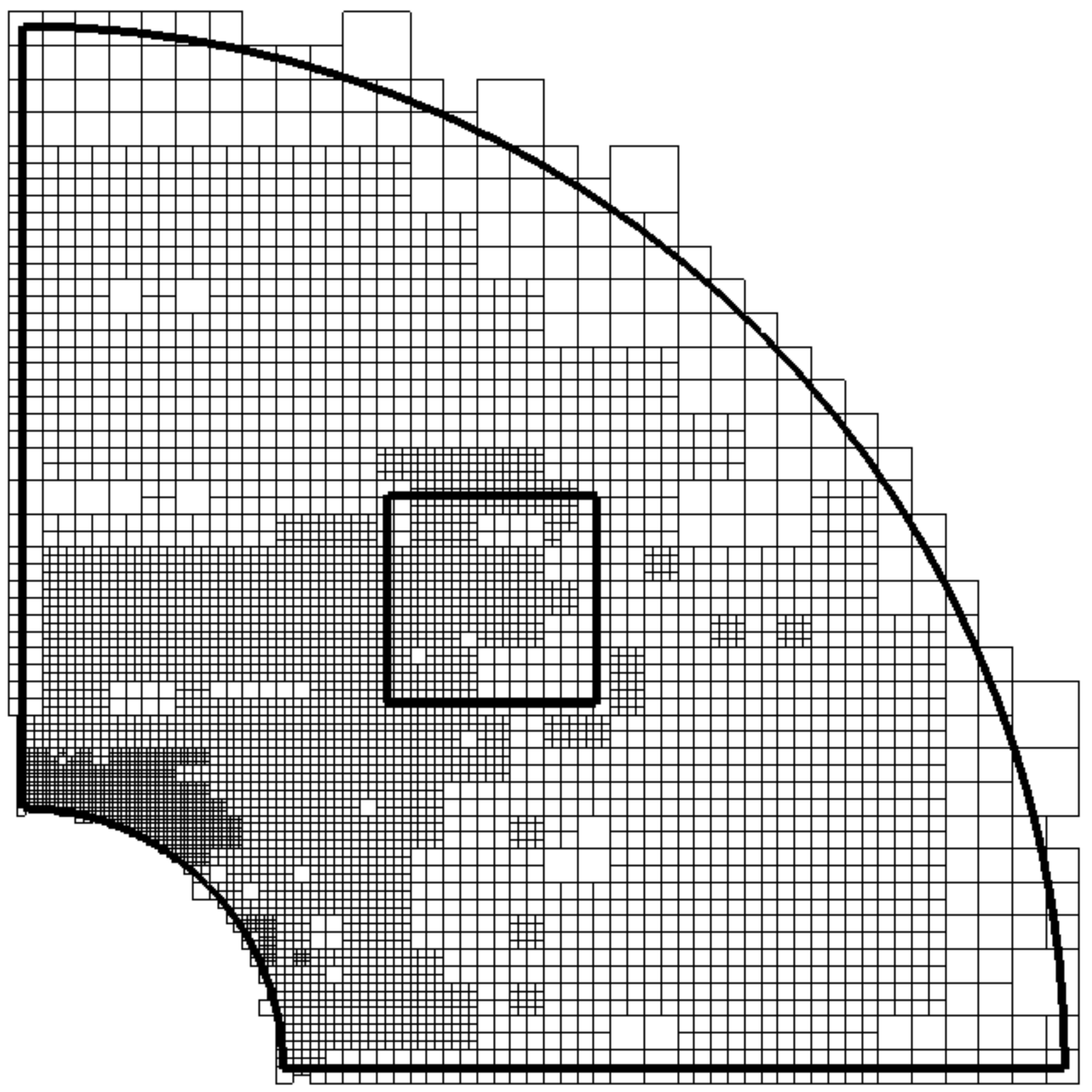} &
\includegraphics[width=0.3\textwidth]{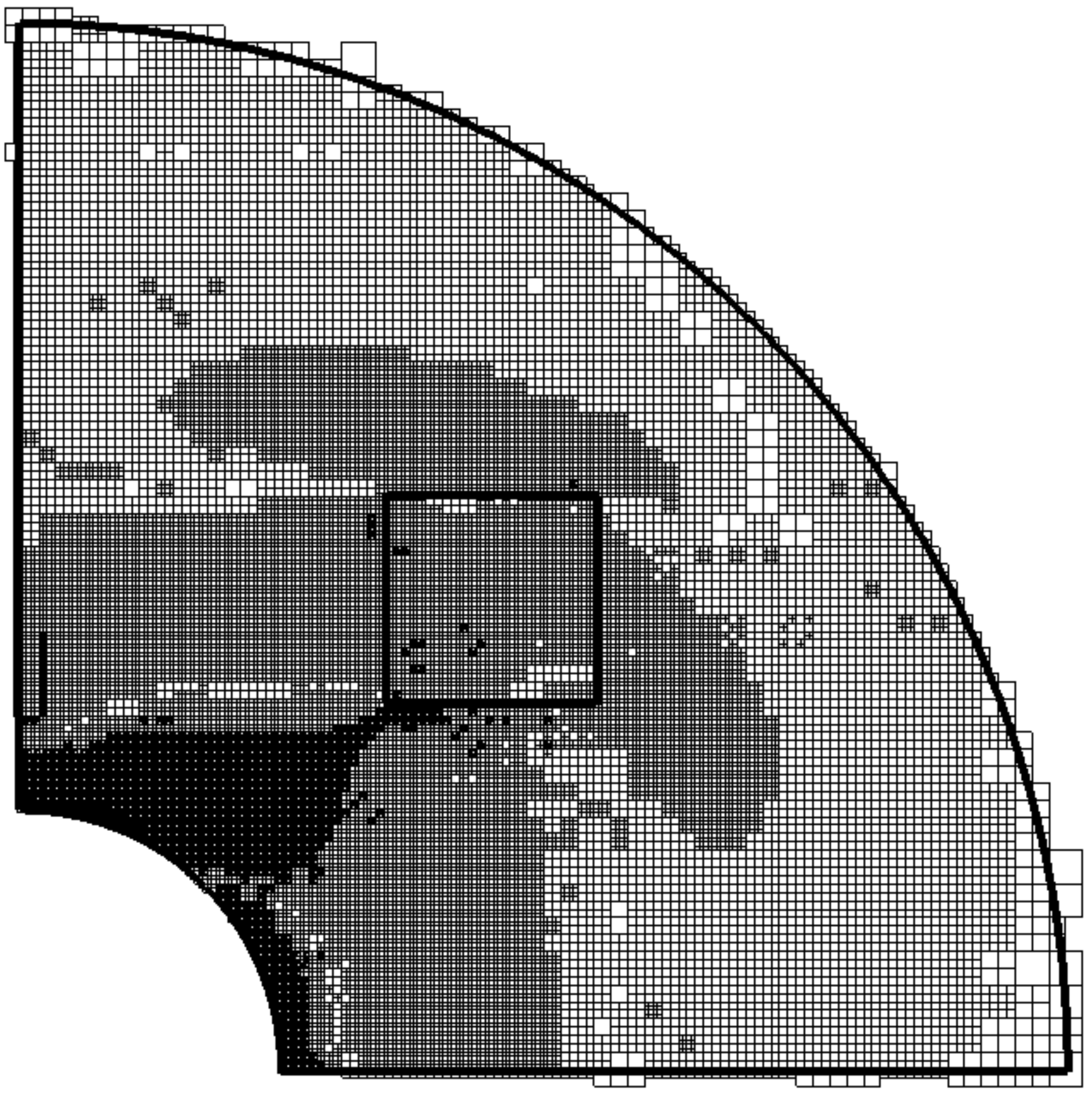}\\
{c) Mesh 3} & {d)Mesh 4}
\end{tabular}
\caption{Problem 1.b. Sequence of \textit{h}-adaptive refined meshes.}%
\label{fig:Pipe_mD_Mesh}%
\end{figure}

The exact value of the QoI can be computed for this problem and is $\bar{u}_x = 0.002238239291713$. Figure \ref{fig:Pipe_mD_Mesh} shows the first four meshes used in the \textit{h}-adaptive refinement process guided by the error estimate for the QoI.  

\begin{figure}[htb!]
      \centering
      \begin{tikzpicture}
    \begin{semilogxaxis}[       
    width=.48\textwidth,
    title={$\theta$},
    xlabel={DOF},
    legend style={cells={anchor=west}, font=\small},
    legend style={at={(0.98,0.98)},anchor=north east},
    cycle list name=ageplot, table/col sep=comma
    ]
   \draw (axis cs: 0,1)--(axis cs: 1e8,1); 
    \addplot table[x=dof,y= theta]{Datos/CYLMeanDisp.csv};
    \addplot table[x=dof,y= E2]{Datos/CYLMeanDisp.csv};
    \addplot table[x=dof,y= E3]{Datos/CYLMeanDisp.csv};
    \addplot table[x=dof,y= E4]{Datos/CYLMeanDisp.csv};
    \legend{{$\mathcal{E}_1$},{$\mathcal{E}_2$},{$\mathcal{E}_3$},{$\mathcal{E}_4$}}
    \end{semilogxaxis}
      \end{tikzpicture}  
\caption{Problem 1.b. Evolution of the effectivity index $\theta$ for the error estimates in  (\ref{Eq:errorEstQoI}, \ref{Eq:errorEstimates2}-\ref{Eq:errorEstimates4}) obtained with the SPR-CX technique.}%
\label{fig:CylThetaDisp}
\end{figure}
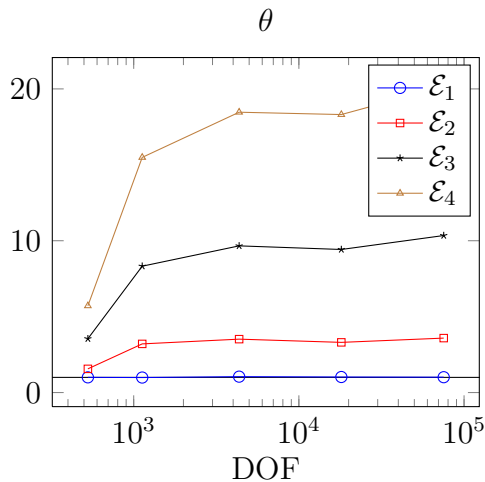

Figure \ref{fig:CylEtaDisp} shows the  relative error (in percentage) for the error estimates in (\ref{Eq:errorEstQoI},\ref{Eq:errorEstimates2}-\ref{Eq:errorEstimates4}) evaluated using the proposed recovery technique. The relative error for $Q(\vm{e})$ is also shown for comparison. For all the curves the relative error decreases monotonically when increasing the number of DOF, indicating that the \emph{h}-adaptive process has a stable convergence. The most accurate estimation is obtained for the estimate $\mathcal{E}_{1}$, which practically coincides with the exact relative error. The other estimates tend to overestimate the exact error, although strictly speaking they do not have bounding properties, as discussed in Section~\ref{sec:errorQoI}.
Figure \ref{fig:CylThetaDisp} represents the evolution of the effectivity index as we increase the number of DOF. The curves are in consonance with the results in Figure \ref{fig:CylEtaDisp} and show good values for the effectivity index, especially for $\mathcal{E}_{1}$ with values within $[0.9898, 1.0481]$. 

\begin{figure}[htb!]
    \centering
    \begin{tikzpicture}
    \begin{loglogaxis}[
    width=0.48\textwidth,
    title={$\eta^Q(\%)$},
    xlabel={DOF},
    legend style={cells={anchor=west}, font=\small},
    legend style={at={(1.02,0.98)},anchor=north west},
    cycle list name=ageplot, table/col sep=comma
    ]
    \draw (axis cs: 0,1)--(axis cs: 1e8,1); 
    \addplot table[x=dof,y= etaes]{Datos/CYLMeanDisp.csv};
    \addplot table[x=dof,y= Eta2]{Datos/CYLMeanDisp.csv};
    \addplot table[x=dof,y= Eta3]{Datos/CYLMeanDisp.csv};
    \addplot table[x=dof,y= Eta4]{Datos/CYLMeanDisp.csv};
    \addplot table[x=dof,y= eta]{Datos/CYLMeanDisp.csv};
    \legend{{$\mathcal{E}_1$},{$\mathcal{E}_2$},{$\mathcal{E}_3$},{$\mathcal{E}_4$},{$Q(\vm{e})$}}
    \end{loglogaxis}
   \end{tikzpicture} 
\caption{Problem 1.b. Evolution of the relative error $\eta^Q$ obtained with the SPR-CX technique, considering the error estimates in  (\ref{Eq:errorEstQoI}, \ref{Eq:errorEstimates2}-\ref{Eq:errorEstimates4}) and the exact error $Q(\vm{e})$.} 
\label{fig:CylEtaDisp}
\end{figure}
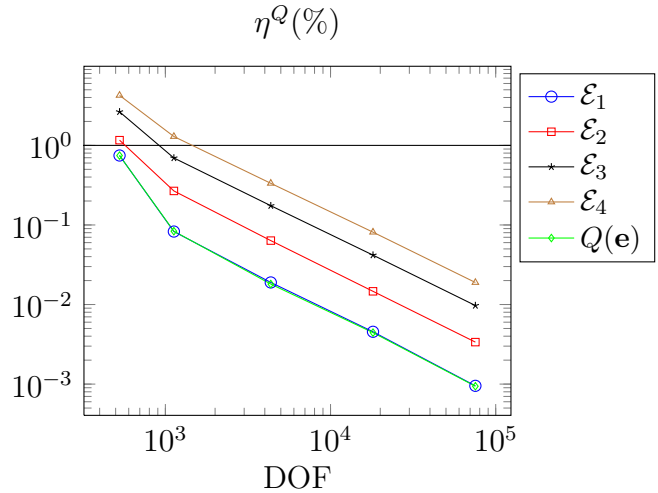

Table \ref{tab:MoI_mD_ConvergQ4} shows the estimated error in the QoI, $\mathcal{E}_1$, the exact error, $Q(\vm{e})$,  the effectivity in the quantity of interest, $\theta_{QoI}$ and the effectivity of the error estimator, $\theta$, for the sharp error estimate $\mathcal{E}_1^{\rm SPR-CX} = \mathcal{E}_1$. Comparing $Q(\vm{e})$ and $\mathcal{E}_1$ we can notice that both values decrease as we refine and that the estimate $\mathcal{E}_1^{\rm SPR-CX}$ gives a good approximation of the exact error. The effectivity of the error estimator $\theta$ converges and is very close to the optimal value $\theta=1$ (with $\theta = 1.0123$ for 75,323 DOF). As expected from these results, the effectivity $\theta_{QoI}$ is very accurate as well, with $\theta_{QoI} = 1.00000012$ for  75,323 DOF. 

\begin{table}
\centering\caption{Problem 1.b. Error estimate $\mathcal{E}_1^{\rm SPR-CX}$ and its effectivities.}
\begin{filecontents*}{Datos/CYLMeanDisp.csv}
dof,Qu,Norm,Qees,MOI,Qe,etaes,eta,theta,Efect0,thetaQoI,E2,E3,E4,Eta2,Eta3,Eta4,GOA_UB,GOA_LB
5.2800000000000E+02,2.2382392917412E-03,2.2215739465489E-03,1.6704293743305E-05,2.2382782402922E-03,1.6665345192324E-05,7.4631402482038E-01,7.4457388241808E-01,1.0023370983638E+00,9.9255426117582E-01,1.0000174014240E+00,1.5622164755740E+00,3.5472715381642E+00,5.7133253045710E+00,1.1631855863956E+00,2.6412057411621E+00,4.2539928035419E+00,3.3578312014674E+00,-2.3554941031036E+00
1.1260000000000E+03,2.2382392917412E-03,2.2363706567901E-03,1.8495994798972E-06,2.2382202562700E-03,1.8686349510496E-06,8.2636360049703E-02,8.3486826361443E-02,9.8981316755222E-01,9.9916513173639E-01,9.9999149533688E-01,3.2062167868710E+00,8.3233023388332E+00,1.5489192284955E+01,2.6767686416264E-01,6.9488609711596E-01,1.2931435067730E+00,8.2395027262535E+00,-7.2496895587012E+00
4.3480000000000E+03,2.2382392917412E-03,2.2378346842623E-03,4.2408593897234E-07,2.2382587702013E-03,4.0460747887652E-07,1.8947301145913E-02,1.8077042985058E-02,1.0481416214794E+00,9.9981922957015E-01,1.0000087025816E+00,3.5181508068135E+00,9.6619781061468E+00,1.8464541415918E+01,6.3597763362686E-02,1.7465999354551E-01,3.3378430887494E-01,9.7563415186987E+00,-8.7081998972193E+00
1.8076000000000E+04,2.2382392917412E-03,2.2381400667799E-03,1.0170695089463E-07,2.2382417737308E-03,9.9224961315567E-08,4.5440606493646E-03,4.4331703800257E-03,1.0250137621235E+00,9.9995566829620E-01,1.0000011089027E+00,3.3060248551546E+00,9.4241090554299E+00,1.8310573367909E+01,1.4656171463500E-02,4.1778681122664E-02,8.1173891495902E-02,9.6677935650163E+00,-8.6427798028928E+00
7.5328000000000E+04,2.2382392917412E-03,2.2382182684281E-03,2.1281959640380E-08,2.2382395503877E-03,2.1023313076084E-08,9.5083486912717E-04,9.3927906429208E-04,1.0123028451015E+00,9.9999060720936E-01,1.0000001155580E+00,3.5869751434999E+00,1.0342919590429E+01,2.0149363113523E+01,3.3691706564256E-03,9.7148878349469E-03,1.8925874931351E-02,1.0580832979312E+01,-9.5685301342105E+00
\end{filecontents*}
\pgfplotstabletypeset[col sep= comma, columns={dof,Qees,Qe,theta,thetaQoI},
 columns/dof/.style   ={column name=dof, column type=r, int detect},
 columns/theta/.style ={column name=$\theta$,fixed, zerofill,precision=8},
 columns/thetaQoI/.style ={column name=$\theta_{QoI}$,fixed, zerofill,precision=8},
 columns/Qees/.style  ={column name={$\mathcal{E}_1^{\rm SPR-CX}$} , zerofill,precision=6},
 columns/Qe/.style    ={column name=$Q(\vm{e})$ , zerofill,precision=6,},
]
{Datos/CYLMeanDisp.csv} 
\label{tab:MoI_mD_ConvergQ4}
\end{table}

\begin{figure}[htb!]
      \centering
      \begin{tikzpicture}
    \begin{semilogxaxis}[
    width=0.48\textwidth,
    title={$\theta$},
    ymin=0.8, ymax=1.3,
    xlabel={DOF},
    legend style={cells={anchor=west}, font=\small},
    legend style={at={(0.98,0.98)},anchor=north east},
    cycle list name=ageplot, table/col sep=comma
    ]a
    \draw (axis cs: 0,1)--(axis cs: 1e8,1); 
    \addplot table[x=dof,y= theta]{Datos/CYLMeanDisp.csv};
    \addplot table[x=dof,y= theta]{Datos/CYLMeanDispSPR.csv};
    \legend{{$\mathcal{E}_1^{\rm SPR-CX}$},{$\mathcal{E}_1^{\rm SPR}$}}
    \end{semilogxaxis}
      \end{tikzpicture}  
\caption{Problem 1.b. Evolution of the effectivity index $\theta$ considering equilibrated $\mathcal{E}_1^{\rm SPR-CX}$ and non-equilibrated recovery, $\mathcal{E}_1^{\rm SPR}$.}  %
\label{fig:CylThetaDispNEq}
\end{figure}
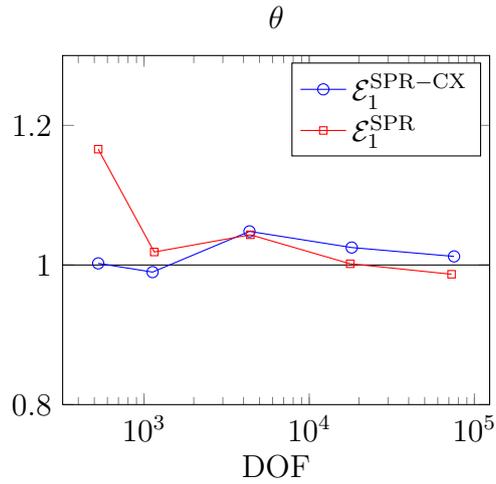

If in (\ref{Eq:errorEstQoI}) we consider the case of the non--equilibrated superconvergent patch recovery procedure, resembling the averaging error estimators presented in \cite{ruterstein2006}, we obtain the results shown in Figure \ref{fig:CylThetaDispNEq}. This figure shows that, in this case, the effectivity of the error estimator provided by the standard SPR technique is similar to the effectivity obtained with the SPR-CX technique here proposed, although the latter results in more accurate values for the coarsest mesh ($\theta = 1.1659$ for the SPR and $\theta = 1.0023$ for the SPR-CX, considering 528 DOF). In any case, we should recall that the local behaviour with the SPR-CX is generally better than with the SPR.

\subsubsection{Problem 1.c.: Mean stress \texorpdfstring{$\bar{\sigma}_x$}{sigmax}  in \texorpdfstring{$\Omega_I$}{Omegai}}

\begin{figure}[htb!]
\centering
\begin{tabular}{c c}
\includegraphics[width=0.3\textwidth]{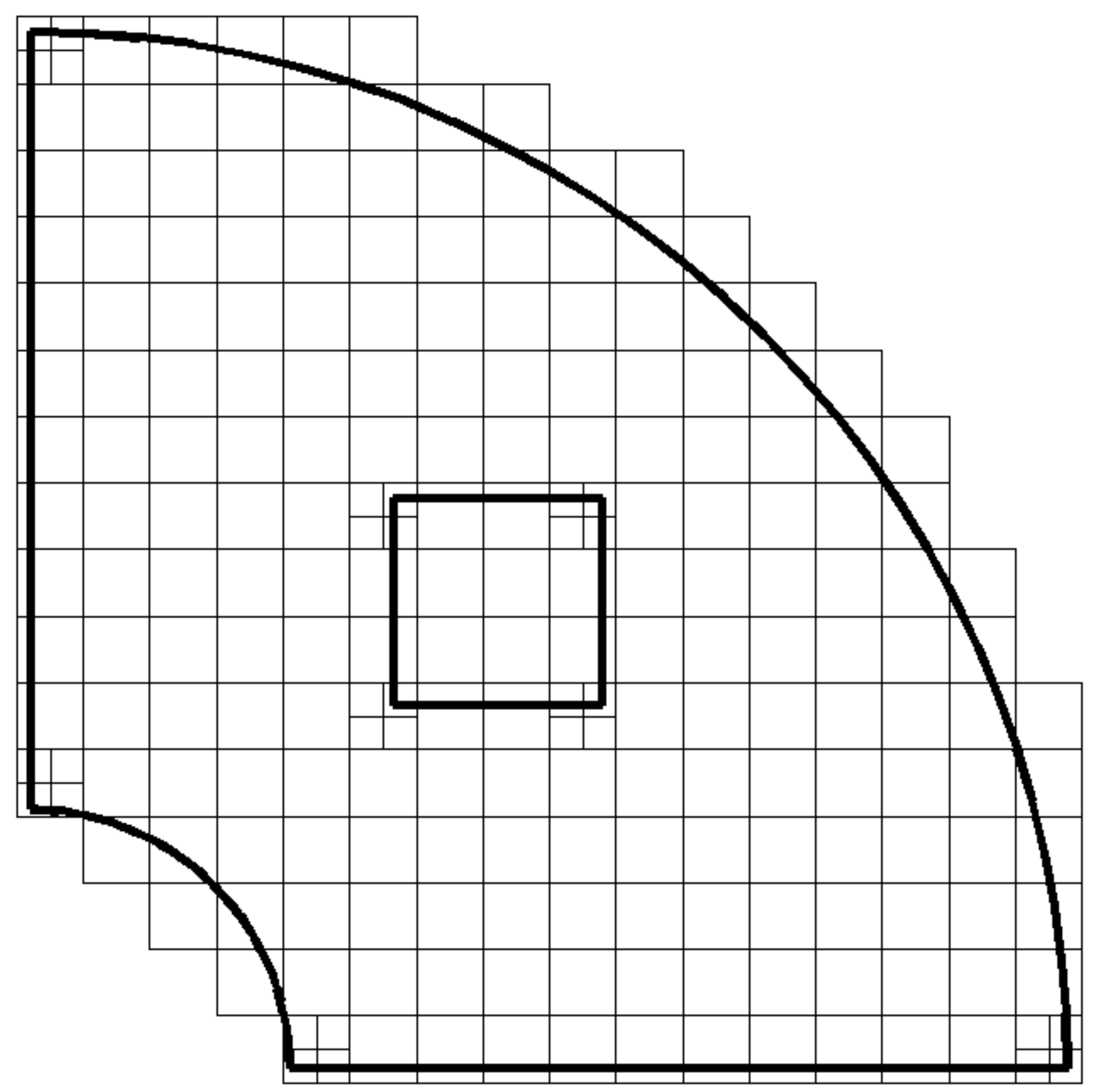}
 &
\includegraphics[width=0.3\textwidth]{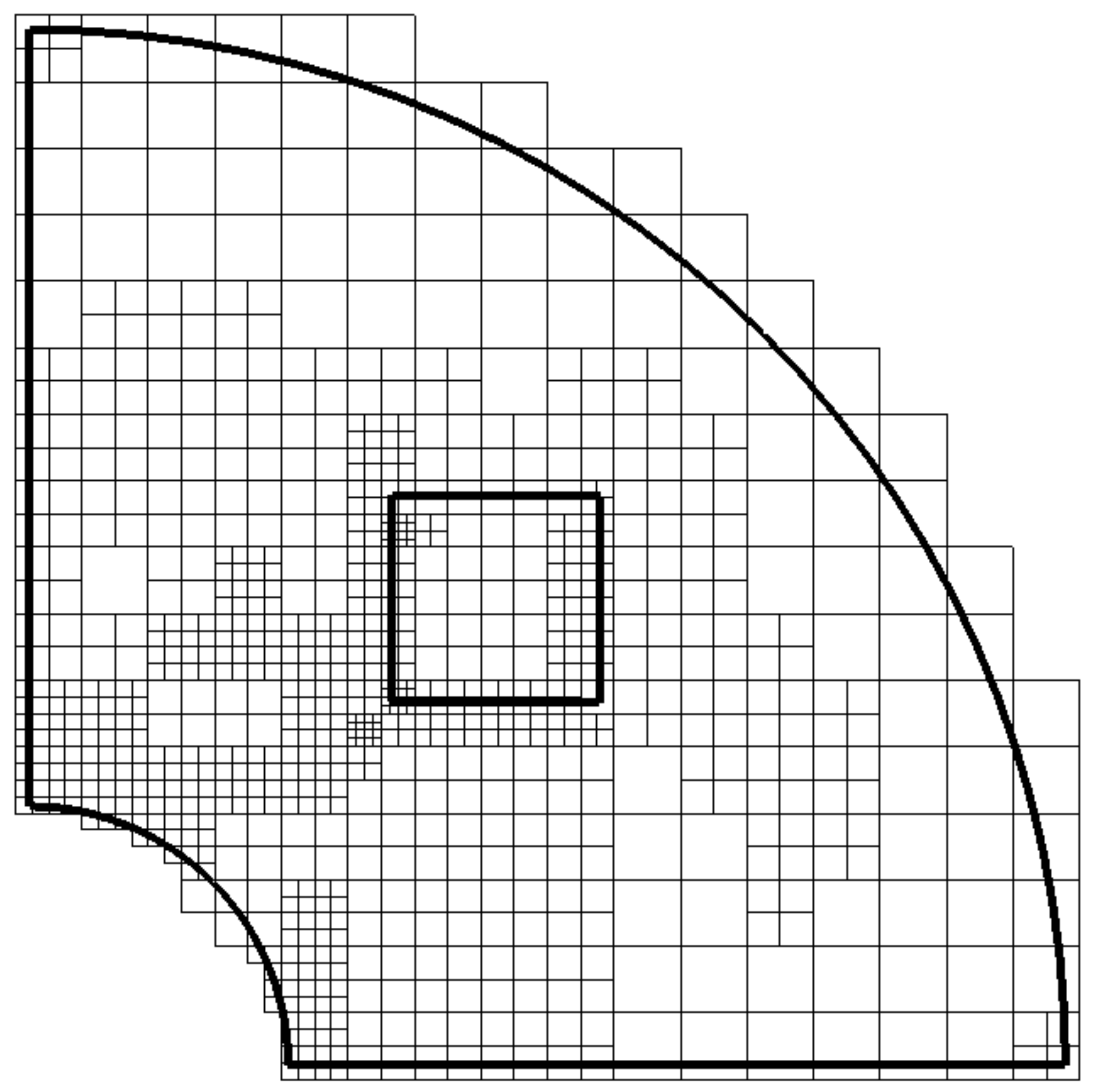}
\\
{a) Mesh 1} & {b) Mesh 2}
\\
\includegraphics[width=0.3\textwidth]{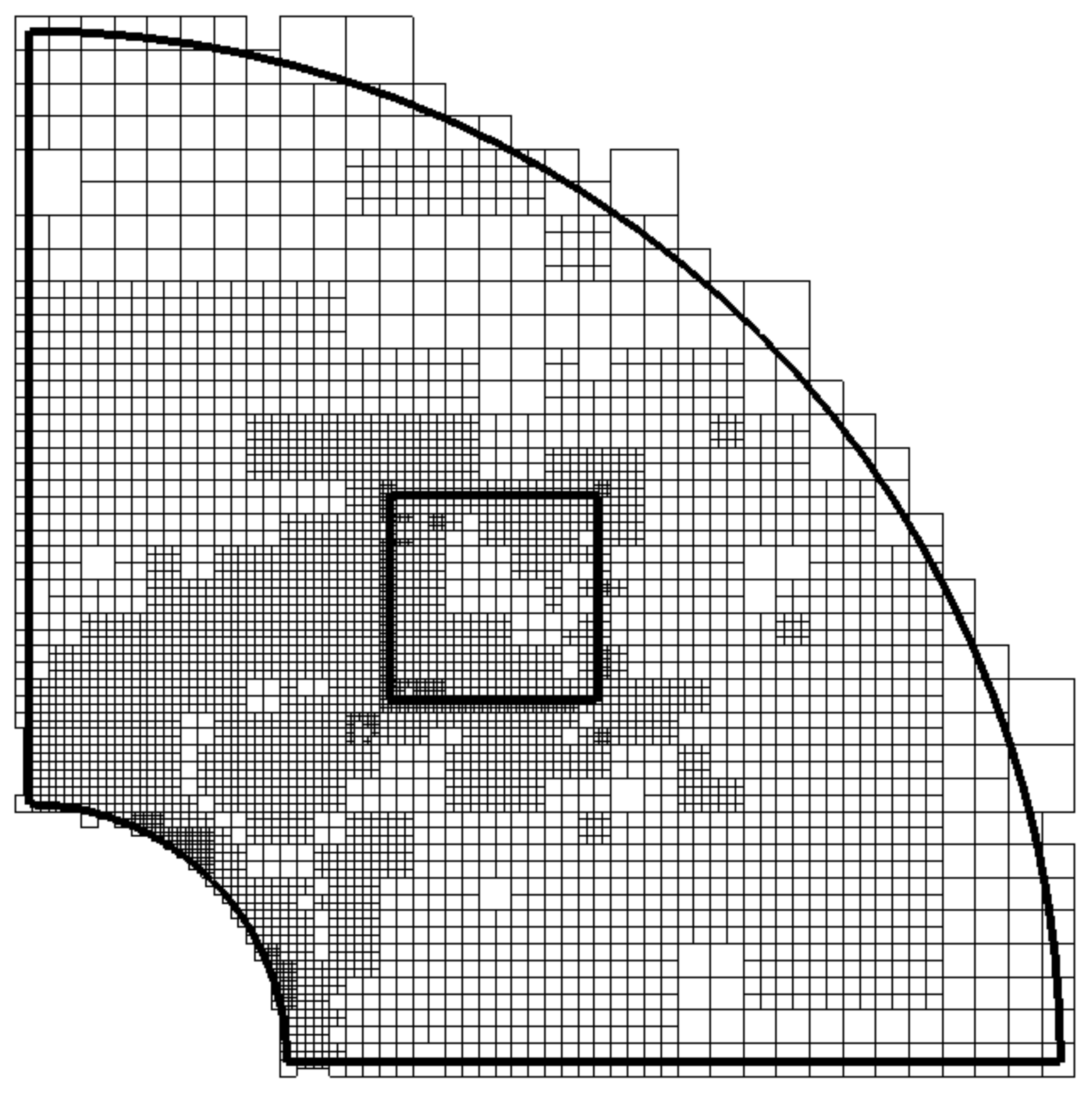}
 &
\includegraphics[width=0.3\textwidth]{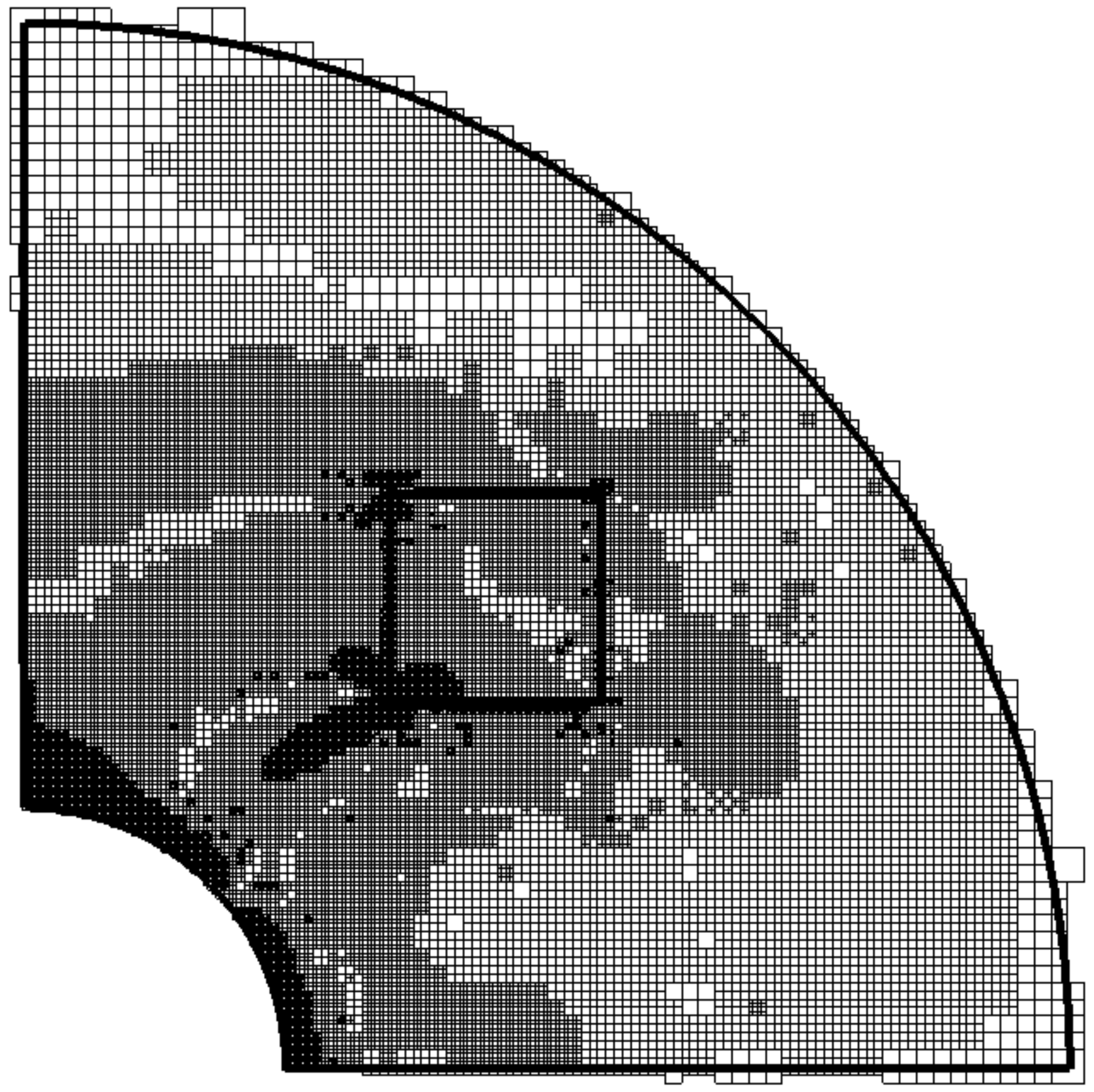}
\\
{c) Mesh 3} & {d) Mesh 4}
\end{tabular}\\
\caption{Problem 1.c. Sequence of meshes for the mean stress in a domain of interest.}%
\label{fig:Pipe_mS_Mesh}%
\end{figure}

Consider as QoI the mean stress value $\bar{\sigma}_x$ given in (\ref{eq:MoI_mStress}). Figure \ref{fig:Pipe_mS_Mesh} shows the first four meshes of bilinear elements used in the refinement process guided by the error estimated for this QoI.

\begin{figure}[htb!] 
    \centering
    \begin{tikzpicture}
    \begin{loglogaxis}[
    width=0.48\textwidth,
    title={$\eta^Q (\%)$},
    xlabel={DOF},
    legend style={cells={anchor=west}, font=\small},
    legend style={at={(1.02,0.98)},anchor=north west},
    cycle list name=ageplot, table/col sep=comma
    ]
    \draw (axis cs: 0,1)--(axis cs: 1e8,1); 
    \addplot table[x=dof,y= etaes]{Datos/CYLMeanStress.csv};
    \addplot table[x=dof,y= Eta2]{Datos/CYLMeanStress.csv};
    \addplot table[x=dof,y= Eta3]{Datos/CYLMeanStress.csv};
    \addplot table[x=dof,y= Eta4]{Datos/CYLMeanStress.csv};
    \addplot table[x=dof,y= eta]{Datos/CYLMeanStress.csv};
    \legend{{$\mathcal{E}_1$},{$\mathcal{E}_2$},{$\mathcal{E}_3$},{$\mathcal{E}_4$},{$Q(\vm{e})$}}
    \end{loglogaxis}
   \end{tikzpicture} 
\caption{Problem 1.c. Evolution of the relative error $\eta^Q$ considering the error estimates in  (\ref{Eq:errorEstQoI}, \ref{Eq:errorEstimates2}-\ref{Eq:errorEstimates4}) and the exact error $Q(\vm{e})$, obtained with the SPR-CX technique.} 
\label{fig:CylEtaStress}
\end{figure}
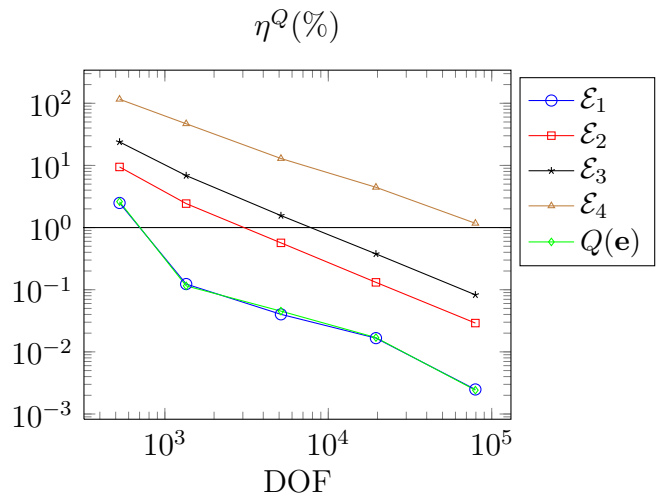

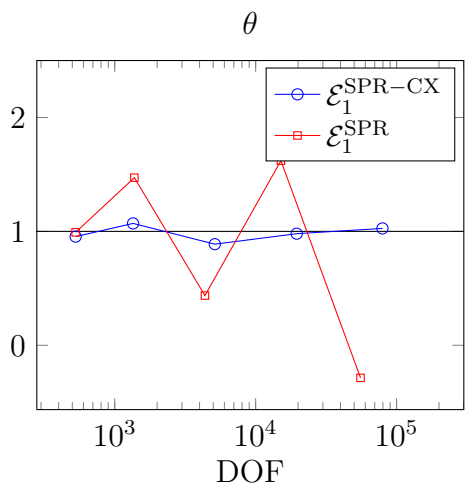
\begin{figure}[htb!]
      \centering
      \begin{tikzpicture}
    \begin{semilogxaxis}[
    width=0.48\textwidth,
    title={$\theta$},
    xmin=0, xmax=3e5,
    ymax=2.5,
    xlabel={DOF},
    legend style={cells={anchor=west}, font=\small},
    legend style={at={(0.98,0.98)},anchor=north east},
    cycle list name=ageplot, table/col sep=comma
    ]a
    \draw (axis cs: 0,1)--(axis cs: 1e8,1); 
    \addplot table[x=dof,y= theta]{Datos/CYLMeanStress.csv};
    \addplot table[x=dof,y= theta]{Datos/CYLMeanStressSPR.csv};
    \legend{{$\mathcal{E}_1^{\rm SPR-CX}$},{$\mathcal{E}_1^{\rm SPR}$}}
    \end{semilogxaxis}
      \end{tikzpicture}  
\caption{Problem 1.c. Evolution of the effectivity index $\theta$ considering equilibrated, $\mathcal{E}_1^{\rm SPR-CX}$, and non-equilibrated recovery, $\mathcal{E}_1^{\rm SPR}$.}  %
\label{fig:CylThetaStress}
\end{figure}

The evolution of the relative error for the estimates presented in (\ref{Eq:errorEstQoI}, \ref{Eq:errorEstimates2}-\ref{Eq:errorEstimates4}) using the proposed recovery technique and the exact error is shown in Figure \ref{fig:CylEtaStress}. Similar to the observations done for the previous examples, the most accurate results are obtained when considering the estimate $\mathcal{E}_{1}$. In this case, the other estimates considerably overestimate the true error. Figure \ref{fig:CylThetaStress} shows the effectivity index for $\mathcal{E}_1^{\rm SPR-CX}$, which uses the locally equilibrated SPR-CX recovery technique, together with the effectivity obtained with the non-equilibrated SPR technique ($\mathcal{E}_1^{\rm SPR}$ curve). This graph clearly shows the improvement obtained with the SPR-CX recovery technique, with effectivities very close to 1 ($\theta = 0.9797$ for the SPR-CX whilst $\theta = 1.6209$ for the SPR, considering 19,573 DOF), in contrast with the oscillatory effectivities provided by the SPR technique.

Table \ref{tab:MoI_mS_ConvergQ4} shows the estimate $\mathcal{E}_1^{\rm SPR-CX}$, the exact error $Q(\vm{e})$, the effectivity index for the QoI $\theta_{QoI}$ and the effectivity of the error estimator $\theta$. For this problem the exact value of the QoI is $\bar{\sigma}_x= 0.0\bar{6}$. Table \ref{tab:MoI_mS_ConvergQ4} indicates that the equilibrated recovery procedure (SPR-CX) provides very accurate estimations of the error in the QoI and in the value of the QoI itself, with $\theta = 1.0259$ and $\theta_{QoI} = 1.00000063$ for 79,442 DOF. 

\begin{table}[htb!]
\centering\caption{Problem 1.c. Values for the error estimate $\mathcal{E}_1^{\rm SPR-CX}$ and effectivities.}
\begin{filecontents*}{Datos/CYLMeanStress.csv}
dof,Qu,Norm,Qees,MOI,Qe,etaes,eta,theta,Efect0,thetaQoI,E2,E3,E4,Eta2,Eta3,Eta4,GOA_UB,GOA_LB
5.2800000000000E+02,6.6666666666667E-02,6.8405589320404E-02,1.6590766708840E-03,6.6746512649520E-02,1.7389226537374E-03,2.4886150063260E+00,2.6083839806061E+00,9.5408307397586E-01,1.0260838398061E+00,1.0011976897428E+00,3.6222787178382E+00,9.1536643339835E+00,4.4509132404106E+01,9.4482937808994E+00,2.3876271412607E+01,1.1609690795354E+02,-2.1777524665065E+01,2.2731607739041E+01
1.3500000000000E+03,6.6666666666667E-02,6.6589795538540E-02,8.2255567222104E-05,6.6672051105762E-02,7.6871128126821E-05,1.2338335083316E-01,1.1530669219023E-01,1.0700450120415E+00,9.9884693307810E-01,1.0000807665864E+00,2.1083866712796E+01,5.9539903642368E+01,4.0600521663740E+02,2.4311109292322E+00,6.8653493423265E+00,4.6815118542436E+01,2.0353763082472E+02,-2.0246758581268E+02
5.1310000000000E+03,6.6666666666667E-02,6.6636496804619E-02,2.6782757754212E-05,6.6663279562373E-02,3.0169862047583E-05,4.0174136631318E-02,4.5254793071374E-02,8.8773219154836E-01,9.9954745206929E-01,9.9994919343560E-01,1.2543427195914E+01,3.4461049249495E+01,2.8566931386559E+02,5.6765020215695E-01,1.5595276528083E+00,1.2927905685829E+01,1.4327852302857E+02,-1.4239079083702E+02
1.9573000000000E+04,6.6666666666667E-02,6.6655348939085E-02,1.1088913480955E-05,6.6666437852566E-02,1.1317727581173E-05,1.6633370221432E-02,1.6976591371759E-02,9.7978268176389E-01,9.9983023408628E-01,9.9999656778850E-01,7.7378880373714E+00,2.2214194913553E+01,2.6199637655056E+02,1.3136296329088E-01,3.7712130970000E-01,4.4478054255804E+00,1.3148807961616E+02,-1.3050829693440E+02
7.9442000000000E+04,6.6666666666667E-02,6.6665053412518E-02,1.6551511552805E-06,6.6666708563674E-02,1.6132541482833E-06,2.4827267329207E-03,2.4198812224249E-03,1.0259704938876E+00,9.9997580118778E-01,1.0000006284551E+00,1.2007397011560E+01,3.4198525583214E+01,4.8421031666547E+02,2.9056474558475E-02,8.2756369893439E-02,1.1717314530032E+00,2.4261814357968E+02,-2.4159217308579E+02
\end{filecontents*}
\pgfplotstabletypeset[columns={dof,Qees,Qe,theta,thetaQoI},
 columns/dof/.style   ={column name=dof, column type=r, int detect},
 columns/theta/.style ={column name=$\theta$,fixed, zerofill,precision=8},
 columns/thetaQoI/.style ={column name=$\theta_{QoI}$,fixed, zerofill,precision=8},
 columns/Qees/.style  ={column name={$\mathcal{E}_1^{\rm SPR-CX}$},
    fixed,sci ,dec sep align, zerofill,precision=6,},
 columns/Qe/.style    ={column name=$Q(\vm{e})$,
    fixed,sci ,dec sep align, zerofill,precision=6,},col sep= comma
]
{Datos/CYLMeanStress.csv} 
\label{tab:MoI_mS_ConvergQ4}
\end{table}

\subsection{Problem 2: Mean normal traction \texorpdfstring{$\bar{t}_n$}{tn}  along  \texorpdfstring{$\Gamma_I$}{Gammai}}

\begin{figure}[htb!]
    \centering
    \includegraphics{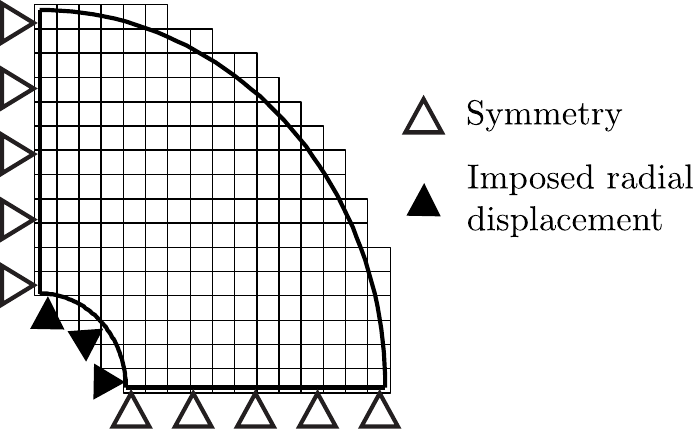}
    \caption{Thick-wall cylinder subjected to an internal pressure.Dirichlet boundary conditions.}
    \label{fig:Pipe_ModelDirichlet}
\end{figure}

In this test case, in the primal problem the displacements are imposed along the inner boundary of the cylinder (see Figure \ref{fig:Pipe_ModelDirichlet}) and the quantity of interest is the mean normal traction along the inner boundary. This is opposite to Problem 1, where tractions were imposed and we were interested in the mean value of radial displacements. In Problem 2, the radial displacements $u_r(a)$ to impose are evaluated from (\ref{Eq:uCylinder}) such that the exact value for the QoI is $\bar{t}_n = 1$.

To formulate the dual problem we consider displacement constraints along $\Gamma_I$ such that
\begin{equation}
u_r(r=a) = \frac{1}{|\Gamma_I|}.  
\end{equation}

\begin{figure}[htb!]
      \centering
      \begin{tikzpicture}
    \begin{semilogxaxis}[
    width=0.48\textwidth,
    title={$\theta$},
    xmin=100, xmax=5000,
    xtick={100,1000},
    xlabel={DOF},
    legend style={cells={anchor=west}, font=\small},
    legend style={at={(0.98,0.02)},anchor=south east},
    cycle list name=ageplot, table/col sep=comma
    ]a
    \draw (axis cs: 0,1)--(axis cs: 1e8,1); 
    \addplot table[x=dof,y= theta]{Datos/CYLMeanSigGamma.csv};
    \addplot table[x=dof,y= theta]{Datos/CYLMeanSigGammaSPR.csv};
    \legend{{$\mathcal{E}_1^{\rm SPR-CX}$},{$\mathcal{E}_1^{\rm SPR}$}}
    \end{semilogxaxis}
      \end{tikzpicture}  
\caption{Problem 2. Evolution of the effectivity index $\theta$ considering locally equilibrated , $\mathcal{E}_1^{\rm SPR-CX}$, and non-equilibrated recovery, $\mathcal{E}_1^{\rm SPR}$.}   
\label{fig:CylThetaSigGamma}
\end{figure}
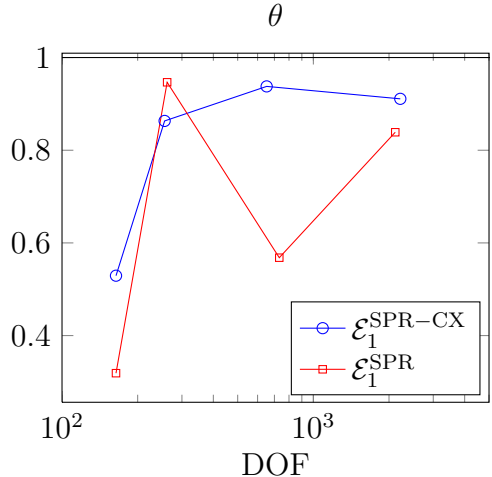

Table \ref{tab:CYLMeanSigGamma} shows the results for the error estimate $\mathcal{E}_1^{\rm SPR-CX}$ using the proposed recovery technique and the exact error $Q(\vm{e})$. The procedure accurately captures the error and yields good effectivities for this QoI, with values of $\theta = 0.9109$ for the last mesh with 2,218 DOF. Figure \ref{fig:CylThetaSigGamma} shows the evolution of the effectivity with the number of DOF for the SPR-CX and the SPR. Results show a better performance of the proposed technique when compared with the SPR which presents an oscillatory behaviour ($\theta = 0.9466$ for mesh 2 with 260 DOF and $\theta = 0.5681$ for mesh 3 with 732 DOF). 

\begin{table}[htb!]
\centering\caption{Problem 2. Values for the error estimate $\mathcal{E}_1^{\rm SPR-CX}$ and effectivities.}
\begin{filecontents*}{Datos/CYLMeanSigGamma.csv}
dof,Qu,Norm,Qees,MOI,Qe,etaes,eta,theta,Efect0,thetaQoI,E2,E3,E4,Eta2,Eta3,Eta4,meanD,sigD
1.6400000000000e+002,1.0000000000000e+000,1.4227986903130e+000,2.2454901512369e-001,1.6473477054367e+000,4.2427618222887e-001,2.2454901512369e-001,4.2427618222887e-001,5.2925199322776e-001,1.4227986903130e+000,1.6473477054367e+000,5.2925199322776e-001,5.2925199322776e-001,5.2925199322776e-001,2.2454901512369e-001,2.2454901512369e-001,2.2454901512369e-001,3.5755920242769e+000,8.1942081659167e+000
2.5600000000000e+002,1.0000000000000e+000,1.0206129999643e+000,1.7940635793597e-002,1.0385536357579e+000,2.0783221124104e-002,1.7940635793597e-002,2.0783221124104e-002,8.6322691205890e-001,1.0206129999643e+000,1.0385536357579e+000,8.6322691205890e-001,8.6322691205890e-001,8.6322691205890e-001,1.7940635793597e-002,1.7940635793597e-002,1.7940635793597e-002,2.5070632352332e-001,3.6137833099858e-001
6.5200000000000e+002,1.0000000000000e+000,1.0041788281517e+000,3.9109830886404e-003,1.0080898112403e+000,4.1705971656676e-003,3.9109830886404e-003,4.1705971656676e-003,9.3775134190270e-001,1.0041788281517e+000,1.0080898112403e+000,9.3775134190270e-001,9.3775134190270e-001,9.3775134190270e-001,3.9109830886404e-003,3.9109830886404e-003,3.9109830886404e-003,1.3579498891913e-001,2.2961471464799e-001
2.2180000000000e+003,1.0000000000000e+000,1.0011068954709e+000,1.0085231941393e-003,1.0021154186651e+000,1.1071656419785e-003,1.0085231941393e-003,1.1071656419785e-003,9.1090542905318e-001,1.0011068954709e+000,1.0021154186651e+000,9.1090542905318e-001,9.1090542905318e-001,9.1090542905318e-001,1.0085231941393e-003,1.0085231941393e-003,1.0085231941393e-003,1.4181307586079e-001,5.5943456251690e-001
\end{filecontents*}
\pgfplotstabletypeset[columns={dof,Qees,Qe,theta,thetaQoI},
 columns/dof/.style   ={column name=dof, column type=r, int detect},
 columns/theta/.style ={column name=$\theta$,fixed, zerofill,precision=8},
 columns/thetaQoI/.style ={column name=$\theta_{QoI}$,fixed, zerofill,precision=8},
 columns/Qees/.style  ={column name={$\mathcal{E}_1^{\rm SPR-CX}$},
    fixed,sci ,dec sep align, zerofill,precision=6,},
 columns/Qe/.style    ={column name=$Q(\vm{e})$,
    fixed,sci ,dec sep align, zerofill,precision=6,},col sep= comma
]
{Datos/CYLMeanSigGamma.csv} 
\label{tab:CYLMeanSigGamma}
\end{table}

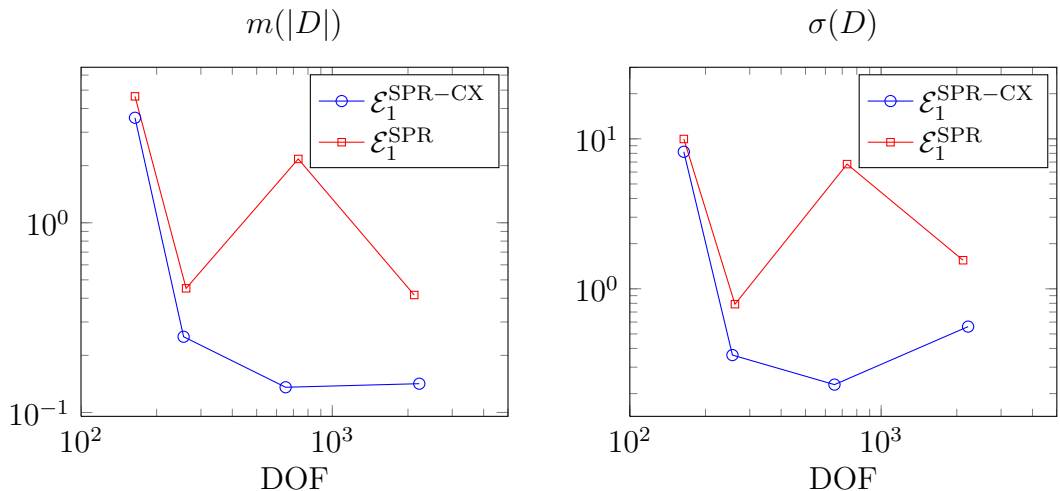
\begin{figure}[htb!]
      \centering
   \begin{minipage}[c]{0.48\textwidth}
       \centering
      \begin{tikzpicture}
    \begin{loglogaxis}[
    width=1\textwidth,
    title={$m(|D|)$},
    xmin=100, xmax=5000,
    xtick={100,1000},
    xlabel={DOF},
    legend style={cells={anchor=west}, font=\small},
    legend style={at={(0.98,0.98)},anchor=north east},
    cycle list name=ageplot, table/col sep=comma
    ]a
    \addplot table[x=dof,y= meanD]{Datos/CYLMeanSigGamma.csv};
    \addplot table[x=dof,y= meanD]{Datos/CYLMeanSigGammaSPR.csv};
    \legend{{$\mathcal{E}_1^{\rm SPR-CX}$},{$\mathcal{E}_1^{\rm SPR}$}}
    \end{loglogaxis}
      \end{tikzpicture}  
   \end{minipage}
   \begin{minipage}[c]{0.48\textwidth}
       \centering
      \begin{tikzpicture}
    \begin{loglogaxis}[
    width=1\textwidth,
    title={$\sigma(D)$},
    xmin=100, xmax=5000,
    xtick={100,1000},
    ymin=0, ymax=3e1,
    xlabel={DOF},
    legend style={cells={anchor=west}, font=\small},
    legend style={at={(0.98,0.98)},anchor=north east},
    cycle list name=ageplot, table/col sep=comma
    ]a
    \addplot table[x=dof,y= sigD]{Datos/CYLMeanSigGamma.csv};
    \addplot table[x=dof,y= sigD]{Datos/CYLMeanSigGammaSPR.csv};
    \legend{{$\mathcal{E}_1^{\rm SPR-CX}$},{$\mathcal{E}_1^{\rm SPR}$}}
    \end{loglogaxis}
      \end{tikzpicture}  
   \end{minipage}
\caption{Problem 2. Evolution of the mean absolute value $m(|D|)$ and standard deviation $\sigma(D)$ of the local effectivity considering locally equilibrated, $ \mathcal{E}_1^{\rm SPR-CX}$, and non-equilibrated recovery, $\mathcal{E}_1^{\rm SPR}$.}   
\label{fig:CylmeanDSigGamma}
\end{figure}

Figure \ref{fig:CylmeanDSigGamma} represents the evolution of the mean absolute  value $m(|D|)$ and standard deviation $\sigma(D)$ of the local effectivity. Again, for this example, the SPR-CX gives lower values of these parameters than the SPR. The improved local performance of the SPR-CX is particularly useful for the adaptive algorithm. Notice that for the second mesh, although the global effectivity is close to unity for the SPR, $m(|D|)$ and $\sigma(D)$ indicate that the SPR-CX is a superior choice when estimating the true error at the element level. Thus, the apparent satisfactory behaviour of the SPR in mesh 2 is due to error compensations of large values of local effectivities.

\subsection{Problem 3: L-Shape plate}

\begin{figure}[!ht]
    \centering
    \includegraphics{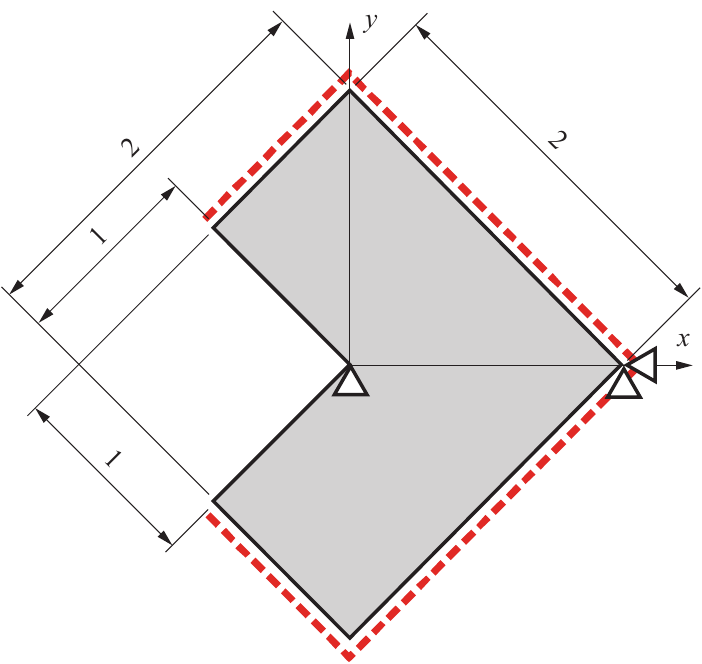}
    \caption{Problem 3. L-shaped domain.}
    \label{fig:LShape}
\end{figure}

Let us consider the singular problem of a finite portion of an infinite domain with a reentrant corner. The model is loaded on the boundary with the tractions corresponding to the first terms of the asymptotic expansion that describes the exact solution under mixed mode loading conditions around the singular vertex, see Figure~\ref{fig:LShape}. The exact values of boundary tractions on the boundaries represented by discontinuous thick lines were imposed in the FE analyses. 

The exact displacement and stress fields for this singular elasticity problem can be found in \cite{szabobabuska1991}. Exact values of the generalised stress intensity factors (GSIF) \cite{szabobabuska1991} under mixed mode were taken as $K_{\rm I}=1$  and $ K_{\rm II}=1$. The material parameters are Young's modulus $E~=~1000$, and Poisson's ratio $\nu~=~0.3$. As the analytical solution of this problem is singular at the reentrant corner of the plate, for the recovery of the dual and primal fields we apply the \emph{singular}+\emph{smooth} decomposition of the stresses as explained in Section \ref{sec:recovery}. We use 
a domain integral method based on extraction functions to obtain an approximation of the recovered singular part as explained in Section \ref{sec:recovery}. 

\begin{figure}[!htbp] 
    \centering
    \includegraphics[scale=0.9]{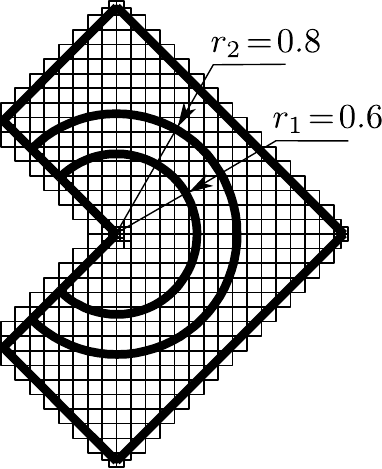}
    \includegraphics[scale=0.9]{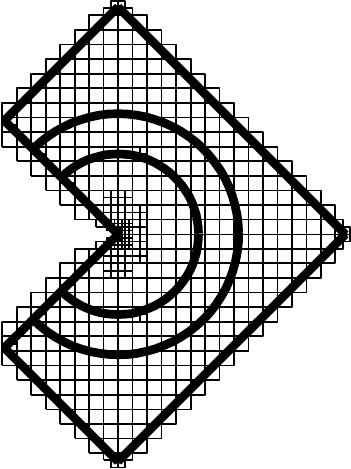}
    \includegraphics[scale=0.9]{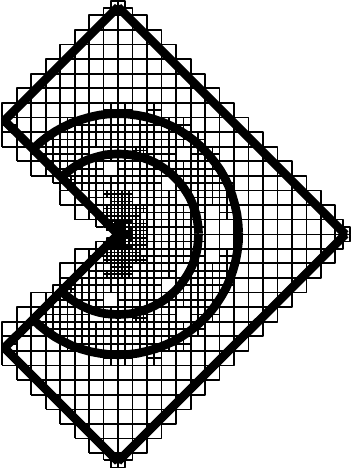}
    \includegraphics[scale=0.9]{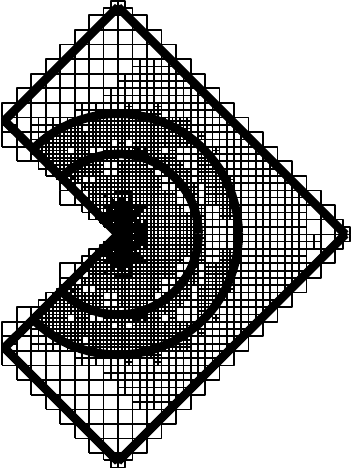}
    \caption{Problem 3. Cartesian meshes with \emph{h}-adaptive refinement.}
    \label{fig:LShapeMesh}
\end{figure}

In this example, we consider the  GSIFs $K_{\rm I}$  and $ K_{\rm II}$ as the quantities of interest. Figure \ref{fig:LShapeMesh} shows the Cartesian meshes used to solve the primal and dual problems when the mesh is \emph{h}-adapted for the evaluation of $K_{\rm I}$. For the dual problem, we use the same Dirichlet conditions as shown in Figure~\ref{fig:LShape} and the set of nodal forces used to extract the QoI in the annular domain $\Omega_I$, defined by a plateau function $q$, shown in Figure~\ref{fig:LShapeMesh}. Function $q$ is defined such that $q=1$ for $r\leq  r_1=0.6$, $q=0$ for $r \geq r_2 =0.8$ and has a smooth transition for $r_1 < r < r_2$ given by a quartic spline $q(s)=1-6s^2+8s^3-3s^4$ with $s=(r-r_1/r_2)/(1-r_1/r_2)$.

In order to impose equilibrium conditions during the recovery of the stresses we use the following approach. For the primal solution, on each patch, we enforce internal equilibrium and compatibility in $\Omega$, and boundary equilibrium all along the Neumann boundary. For the dual problem, we enforce internal equilibrium  using the initial strains and body loads given by (\ref{Eq:IniStrainSIF}) and homogeneous Neumann boundary conditions. Compatibility is enforced within all the domain. 

Table \ref{tab:LSHPM1SPRCX} shows the results for the stress intensity factor $K_{\rm I}$. Similarly to the results for other QoIs, we observe that the proposed technique provides an accurate estimate $\mathcal{E}$ of the exact error $Q(\vm{e})$. The effectivity index $\theta$ is always close to the optimal value $\theta=1$ and for the last mesh with 39,193 DOF we obtain $\theta = 0.9940$. As a result, and in agreement with the previous cases, the effectivity in the QoI is highly accurate, with $\theta_{QoI} = 0.9999$ for the same mesh. Table \ref{tab:LSHPM2SPRCX} shows the same results for the stress intensity factor $K_{\rm II}$. Again, we observe a satisfactory behaviour of the error indicator and very accurate effectivities, both for the error estimate and for the QoI itself ($\theta = 1.0222$ and $\theta_{QoI} = 1.00000053$ for the last mesh).

\begin{table}
\centering\caption{Problem 3. Stress intensity factor $K_{\rm I}$ as QoI using the initial strains and body loads in \eqref{Eq:IniStrainSIF}.}
\begin{filecontents*}{Datos/LSHPM1SPRCX.csv}
dof,theta_p,theta,Kp_1,Kp_2,Kp_(FU),Kd_1,Kd_2,Qees,Qe,thetaQoI,theta2,theta3,theta4,Eta2,Eta3,Eta4 ,etaes,eta
1.1010000000000e+003,9.2247250012288e-001,7.6335570420301e-001,9.9043324898809e-001,9.9820228433630e-001,9.9043324898810e-001,1.1286265863987e+002,4.8865684135525e-011,7.3028339556238e-003,9.5667510119000e-003,9.9773608294372e-001,1.6642694066047e+000,6.2291468045034e+000,1.8542036456339e+001,1.5921651029710e-002,5.9592696495256e-002,1.7738704603137e-001,0.00731940448026846,0.0095667510119
2.9270000000000e+003,9.7363971455278e-001,9.3631773936489e-001,9.9779365690239e-001,9.9956575419468e-001,9.9779365690238e-001,1.1441866794249e+002,6.3831061098509e-003,2.0658381814233e-003,2.2063430976161e-003,9.9985949508381e-001,1.8292177047638e+000,7.0385336960280e+000,2.7237470163376e+001,4.0358818569427e-003,1.5529420237570e-002,6.0095204291490e-002,0.0020661284826326,0.0022063430976161
1.0399000000000e+004,9.9100307888021e-001,9.8377706728170e-001,9.9947198384470e-001,9.9990052888356e-001,9.9947198384468e-001,1.1493505940818e+002,7.6387232168226e-003,5.1945018476238e-004,5.2801615532438e-004,9.9999143402944e-001,1.9473153492279e+000,5.9695477072760e+000,2.5344227312300e+001,1.0282139639035e-003,3.1520176294214e-003,1.3382161465108e-002,0.000519454634395484,0.00052801615532438
3.9193000000000e+004,9.9641107012433e-001,9.9409058406738e-001,9.9986921380313e-001,9.9997612859093e-001,9.9986921380317e-001,1.1506628194008e+002,2.3657888543260e-003,1.3001332679528e-004,1.3078619683060e-004,9.9999922712996e-001,1.9460406072335e+000,4.9458207485769e+000,2.2731734660713e+001,2.5451524989798e-004,6.4684508591224e-004,2.9729971236370e-003,0.000130013427278762,0.0001307861968306
\end{filecontents*}
\pgfplotstabletypeset[columns={dof,Qees,Qe,theta,thetaQoI},
 columns/dof/.style   ={column name=dof, column type=r, int detect},
 columns/theta/.style ={column name=$\theta$,fixed, zerofill,precision=8},
 columns/thetaQoI/.style ={column name=$\theta_{QoI}$,fixed, zerofill,precision=8},
 columns/Qees/.style  ={column name={$\mathcal{E}_1^{\rm SPR-CX}$},
    fixed,sci ,dec sep align, zerofill,precision=6,},
 columns/Qe/.style    ={column name=$Q(\vm{e})$,
    fixed,sci ,dec sep align, zerofill,precision=6,},col sep= comma
]
{Datos/LSHPM1SPRCX.csv} 
\label{tab:LSHPM1SPRCX}
\end{table}

\begin{table}
\centering\caption{Problem 3. Stress intensity factor $K_{\rm II}$ as QoI using the initial strains and body loads in \eqref{Eq:IniStrainSIF}.}
\begin{filecontents*}{Datos/LSHPM2SPRCX.csv}
dof,theta_p,theta,Kp_1,Kp_2,Kp_(FU),Kd_1,Kd_2,Qees,Qe,thetaQoI,theta2,theta3,theta4,Eta2,Eta3,Eta4 ,etaes,eta
1.1010000000000e+003,9.2247250012288e-001,1.1300165895373e+000,9.9043324898809e-001,9.9820228433630e-001,9.9820228433630e-001,1.1066379396100e-010,7.8126030040686e+002,2.0314485232538e-003,1.7977156637016e-003,1.0002337328596e+000,3.7533197592705e+001,2.0820689415135e+002,6.8157207427512e+002,6.7474017221212e-002,3.7429679490654e-001,1.2252727938660e+000,0.00203097381793566,0.0017977156637016
2.9270000000000e+003,9.7363971455278e-001,1.2143971326978e+000,9.9779365690239e-001,9.9956575419468e-001,9.9956575419468e-001,6.3665311217360e-003,7.9069422395963e+002,5.2734686086980e-004,4.3424580532259e-004,1.0000931010555e+000,3.5706730796939e+001,2.0821604656644e+002,9.0232920149965e+002,1.5505498070353e-002,9.0416944822328e-002,3.9183267077130e-001,0.000527297768890927,0.00043424580532259
1.0399000000000e+004,9.9100307888021e-001,1.0188199343063e+000,9.9947198384470e-001,9.9990052888356e-001,9.9990052888354e-001,7.6276506376094e-003,7.9327731567298e+002,1.0134315633345e-004,9.9471116456362e-005,1.0000018720399e+000,3.2321213435750e+001,1.5530555896476e+002,7.4990304552277e+002,3.2150271856785e-003,1.5448417342104e-002,7.4593693172176e-002,0.000101342966615376,0.000099471116456362
3.9193000000000e+004,9.9641107012433e-001,1.0221967537771e+000,9.9986921380313e-001,9.9997612859093e-001,9.9997612859093e-001,2.2485130999010e-003,7.9396762001310e+002,2.4401276857178e-005,2.3871409067788e-005,1.0000005298678e+000,2.9692209481747e+001,1.1276767189377e+002,6.2540519044119e+002,7.0879487866524e-004,2.6919232253983e-003,1.4929303134140e-002,0.0000244012639277343,0.000023871409067788
\end{filecontents*}
\pgfplotstabletypeset[columns={dof,Qees,Qe,theta,thetaQoI},
 columns/dof/.style   ={column name=dof, column type=r, int detect},
 columns/theta/.style ={column name=$\theta$,fixed, zerofill,precision=8},
 columns/thetaQoI/.style ={column name=$\theta_{QoI}$,fixed, zerofill,precision=8},
 columns/Qees/.style  ={column name={$\mathcal{E}_1^{\rm SPR-CX}$},
    fixed,sci ,dec sep align, zerofill,precision=6,},
 columns/Qe/.style    ={column name=$Q(\vm{e})$,
    fixed,sci ,dec sep align, zerofill,precision=6,},col sep= comma
]
{Datos/LSHPM2SPRCX.csv} 
\label{tab:LSHPM2SPRCX}
\end{table}

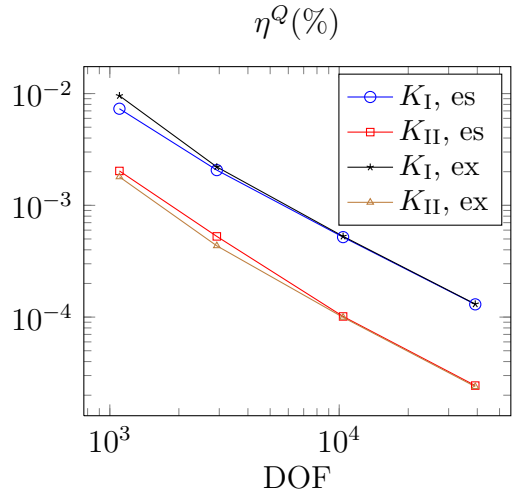
\begin{figure}[htb!]
      \centering
      \begin{tikzpicture}
    \begin{loglogaxis}[
    width=0.48\textwidth,
    title={$\eta^Q(\%)$},
    xlabel={DOF},
    legend style={cells={anchor=west}, font=\small},
    legend style={at={(0.98,0.98)},anchor=north east},
    cycle list name=ageplot, table/col sep=comma
    ]
    \draw (axis cs: 0,1)--(axis cs: 1e8,1); 
    \addplot table[x=dof,y= etaes]{Datos/LSHPM1SPRCX.csv};
    \addplot table[x=dof,y= etaes]{Datos/LSHPM2SPRCX.csv};
    \addplot table[x=dof,y= eta]{Datos/LSHPM1SPRCX.csv};
    \addplot table[x=dof,y= eta]{Datos/LSHPM2SPRCX.csv};   
    \legend{{$K_{\rm I}$, es},{$K_{\rm II}$, es},{$K_{\rm I}$, ex},{$K_{\rm II}$, ex}}
 \end{loglogaxis}
      \end{tikzpicture}  
\caption{Problem 3. Evolution of the relative error $\eta^Q$ for QoIs $K_{\rm I}$ and $K_{\rm II}$, for the exact and estimated error as defined in (\ref{Eq:RelativeErrors}), obtained with the SPR-CX technique.} 
\label{fig:LSHPEta}
\end{figure}

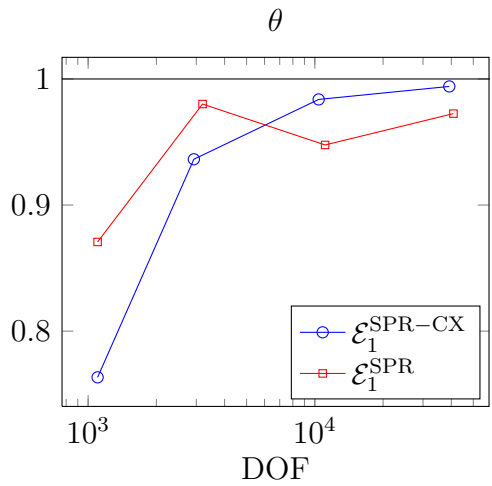
\begin{figure}[htb!]
      \centering
      \begin{tikzpicture}
    \begin{semilogxaxis}[
    width=0.48\textwidth,
    title={$\theta$},
    xlabel={DOF},
    legend style={cells={anchor=west}, font=\small},
    legend style={at={(0.98,0.02)},anchor=south east},
    cycle list name=ageplot, table/col sep=comma
    ]
    \draw (axis cs: 0,1)--(axis cs: 1e8,1); 
    \addplot table[x=dof,y= theta]{Datos/LSHPM1SPRCX.csv};
    \addplot table[x=dof,y= theta]{Datos/LSHPM1SPR.csv};
    \legend{{$\mathcal{E}_1^{\rm SPR-CX}$}, {$\mathcal{E}_1^{\rm SPR}$} }
    \end{semilogxaxis}
      \end{tikzpicture}  
\caption{Problem 3. Evolution of the effectivity index $\theta$ for $K_{\rm I}$} 
\label{fig:LSHPM1_ThetaCompare}
\end{figure}

\begin{figure}[htb!]
      \centering
      \begin{tikzpicture}
    \begin{semilogxaxis}[
    width=0.48\textwidth,
    title={$\theta$},
    xlabel={DOF},
    legend style={cells={anchor=west}, font=\small},
    legend style={at={(0.98,0.98)},anchor=north east},
    cycle list name=ageplot, table/col sep=comma
    ]
    \draw (axis cs: 0,1)--(axis cs: 1e8,1); 
    \addplot table[x=dof,y= theta]{Datos/LSHPM2SPRCX.csv};
    \addplot table[x=dof,y= theta]{Datos/LSHPM2SPR.csv};
    \legend{{$\mathcal{E}_1^{\rm SPR-CX}$}, {$\mathcal{E}_1^{\rm SPR}$} }
    \end{semilogxaxis}
      \end{tikzpicture}  
\caption{Problem 3. Evolution of the effectivity index $\theta$ for $K_{\rm II}$} 
\label{fig:LSHPM2_ThetaCompare}
\end{figure}
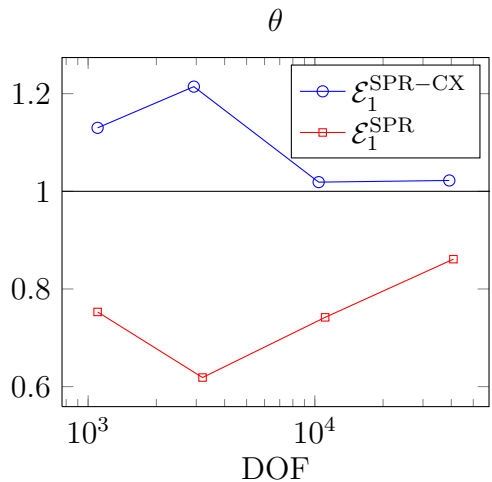

Figure \ref{fig:LSHPEta} shows the evolution of the relative error with respect to the number of DOF for the two QoIs. Figures \ref{fig:LSHPM1_ThetaCompare} and \ref{fig:LSHPM2_ThetaCompare} show the evolution of the effectivity index of the error estimators obtained with the locally equilibrated SPR-CX technique, and with the non-equilibrated recovery as we increase the number of DOF. The results indicate that the proposed methodology accurately evaluates the error in the QoI, giving values of $\theta$ close to 1 and considerably improving the results obtained with the original SPR technique (for the finest mesh considering $K_{\rm II}$ we have $\theta = 0.8610$ for the standard SPR and $\theta = 1.0222$ for the SPR-CX).

\section{Conclusions and future work}

In this paper, we presented an \emph{a posteriori} recovery-based strategy that aims to control the error in quantities of interest. The proposed technique is the first attempt to use a recovery procedure that constructs locally equilibrated stress fields for both the primal and the dual problem.

To recover the solution for the primal problem the formulation enforces equilibrium and compatibility and, for singular problems, relies on a \textit{singular}+\textit{smooth} stress splitting \cite{rodenasgonzalez2008}. In order to enforce equilibrium conditions when evaluating the recovered stress fields for the dual problem, we used an approach that allows us to express the functional which defines a given QoIs in terms of body loads, boundary tractions and initial strains and stresses applied to the dual problem. The proposed technique was tested on different quantities of interest: mean displacements and stresses on a domain of interest, mean displacements and tractions along a boundary and the generalised stress intensity factor for a particular singular problem. 

The methodology we proposed provides accurate global and local  evaluations of the error in the different quantities of interest analysed, improving the results obtained with the original SPR technique, which in some cases can produce satisfactory estimations of the global error due to locally compensating errors. In particular, our approach proves superior to the standard SPR in terms of local error. We also showed that these accurate, local, estimations of the approximation error can be used to drive adaptive procedures.

An extension of the work presented here for extended finite element approximations is now under development. The final goal is to be able to guide the adaptive process in 3D XFEM problems, modelling fatigue crack growth using the stress intensity factors as design parameters for industrial applications of the XFEM, similar to those discussed in the early work of References \cite{bordasmoran2006, bordasconley2007, bordasnguyen2007}  and \cite{wyartcoulon2007,wyartduflot2008}.

\section{Acknowledgements}

This work was supported by the EPSRC grant EP/G042705/1 ``Increased Reliability for Industrially Relevant Automatic Crack Growth Simulation with the eXtended Finite Element Method''. 

St\'ephane Bordas also thanks partial funding for his time provided by the European Research Council Starting Independent Research Grant (ERC Stg grant agreement No. 279578) ``RealTCut Towards real time multiscale simulation of cutting in non-linear materials with applications to surgical simulation and computer guided surgery''. 

This work has received partial support from the research project DPI2010-20542 of the Ministerio de Econom\'ia y Competitividad (Spain). The financial support of the FPU program (AP2008-01086), the funding from Universitat Polit\`{e}cnica de Val\`{e}ncia and Generalitat Valenciana (PROMETEO/2012/023) are also acknowledged. 

All authors also thank the partial support of the Framework Programme 7 Initial Training Network Funding under grant number 289361 ``Integrating Numerical Simulation and Geometric Design Technology."

\bibliographystyle{libstyle}
\bibliography{library}

\clearpage
\end{document}